Dear Colleagues,

My article attached to this email has already been published (information about the journal is given at the top of the first page of the article, and DOI at the bottom of each sheet). But when converting MS word-format, in which the journal demanded to present the material, into the pdf-format, in which the article was published, the journal used a bad translator and did not check the compliance of the resulting file with the original one. In the result, some operations are not reflected or are incorrectly reflected in the pdf-file and this file does not correspond to the source MS word-file. I made my own pdf-file, checked it for compliance with the original MS word-file and asked Editor to replace in the journal electronic version the incorrect pdf-file (not corresponding to the source MS word-file) of the published article with my correct pdf-one (completely corresponding to the source MS word-file). My request was not taken into account. For this reason, I publish in additional open sources the correct pdf-file of my already published article.

Professor A. Rogovoy

# Differentiation of Scalar and Tensor Functions of Tensor Argument

Anatoli A. Rogovoy
*Institute of Continuous Media Mechanics, Russian Academy of Sciences*
*Acad. Korolev Str. 1, 614013 Perm, Russia*
*Corresponding Author: Anatoli A. Rogovoy*

**Abstract:** *In this paper, we analyze the existing rules for constructing derivatives of the scalar and tensor functions of the tensor argument with respect to the tensor argument and the theoretical positions underlying the construction of these rules. We perform a comparative analysis of these rules and the results obtained in the framework of these rules. Considering the existing approaches, we pay due attention to the earliest of them which for some reason is not reflected in later publications on the issue under consideration, and we give to this approach the further development.*

*The rules for constructing the derivatives of scalar and tensor functions of the tensor argument with respect to a tensor and the form of representing these derivatives depend on the accepted scheme of the interaction between the basis vectors in a double scalar product of two tensors. Up to now, three groups of rules and forms of derivative representation have been developed and used. For any of these groups, the derivatives of the scalar functions of the tensor argument with respect to a tensor are similar. The only difference is in derivatives of the tensor functions of the tensor argument.*

*We consider three schemes of the interaction of the basis vectors in a double scalar product of the two tensors, which lie at the root of the rules for constructing derivatives and derivatives themselves and establish the correlation between these schemes. We introduce three isotropic fourth-rank tensors, which play the role of the unit tensor and the unit tensor with the transpose operation for the relevant scheme of the double scalar product by the second-rank tensors. For these schemes we formulate the rules for differentiation of the scalar and tensor functions of the tensor argument with respect to a tensor and for one of them (used in our previous works) we derive specific forms of representation, which reflect the rules for constructing derivatives and expressions for derivatives themselves. In conclusion, we discuss the rules for differentiation of tensor functions of the tensor argument and derivatives obtained in the publications of other authors and examine them for compliance with the construction rules and derivatives obtained in this article.*

## I. Introduction

In solid mechanics, differentiation of the scalar or tensor (second rank) functions of the tensor (also of second rank) argument with respect to a tensor is an often-used procedure. The need of finding derivatives of such functions arises, for example (see [1,2]),

- when deriving from the second law of thermodynamics the constitutive equations, describing the thermomechanical behavior of the material (it is necessary to differentiate the scalar function – Helmholtz free energy with respect to the second-rank tensor, which is any kinematic quantity);
- when constructing the heat equation from the first law of thermodynamics (it is necessary to find time derivative of the scalar function – entropy, depending on the force and kinematic tensors, temperature and parameters, describing changes in the structure of the material under deformation);
- or when developing the relations describing the history of the deformation process step-by-step (relations in terms of increments, the linearization procedure of any scalar or tensor quantity).

However, in the known courses of tensor analysis [3-6] the issues concerning the differentiation of the scalar and tensor functions of the tensor argument with respect to a tensor are usually not considered or poorly elucidated [7]. This is a reason that each time the derivatives of any particular function are found using the approaches of the theory of infinitesimals (increments of the function and the argument are built and passages to the limit are carried out), but not the rules for constructing derivatives, which could have been obtained in the framework of this theory. The rules of the conventional mathematical analysis for finding derivatives of the tensor function of the tensor argument are generally not applicable.

In this paper, we will analyze the existing rules for constructing derivations of the scalar and tensor functions of the tensor argument with respect to the tensor argument and the theoretical positions underlying the construction of these rules. We will perform a comparative analysis of these rules and the results obtained in the framework of these rules. Considering the existing approaches, we will pay due attention to the earliest of them, which for some reason is not reflected in later publications on the issue under consideration, and will give to this approach the further development.

In the literature, there has been a series of works [8-18] that are devoted to the development of the rules for constructing derivatives of the scalar and tensor functions of the tensor argument with respect to a tensor. They also give some examples of how such derivatives are defined. In these works the construction of derivatives is implemented with the use of tensor product operators $\otimes, \boxtimes, \hat{\boxtimes}$ (see, for example, [8-10, 14]), which have the following properties:

$$\underset{\substack{|\,|\\1\,2}}{\mathbf{A}} \otimes \underset{\substack{|\,|\\3\,4}}{\mathbf{B}}, \qquad \underset{\substack{|\,|\\1\,3}}{\mathbf{A}} \boxtimes \underset{\substack{|\,|\\2\,4}}{\mathbf{B}}, \qquad \underset{\substack{|\,|\\1\,4}}{\mathbf{A}} \hat{\boxtimes} \underset{\substack{|\,|\\2\,3}}{\mathbf{B}}, \tag{1.1}$$

where the vertical bars in front of each of the two second-rank tensors $\mathbf{A}$ and $\mathbf{B}$ correspond to the first (left) and second (right) basis vectors, to which these tensors are assigned in the covariant, contravariant or mixed representations, whereas the numbers below them denote the sequence order of these basis vectors in the resulting fourth-rank tensor. Based on the unit isotropic second-rank tensor $\mathbf{g}$ and using the operations $\otimes, \boxtimes$ and $\hat{\boxtimes}$, we also introduce the fourth-rank tensors $\mathbb{I}$ and $\mathbb{T}$, which play the role of the unit tensor and the unit tensor with a transposition operation* in their a double scalar product by the second-rank tensors. The latter is realized according to certain schemes of the interaction between the basis vectors of these fourth-rank tensors and the second-rank tensors (this will be discussed below), which are supposed to comply with the operations in the following relations:

$$\mathbf{A} \cdot\cdot \mathbb{I} = \mathbb{I} \cdot\cdot \mathbf{A} = \mathbf{A}, \qquad \mathbf{A} \cdot\cdot \mathbb{T} = \mathbb{T} \cdot\cdot \mathbf{A} = \mathbf{A}^T. \tag{1.2}$$

The introduced operators and tensors allow us to construct derivatives of the scalar and tensor functions of the tensor argument with respect to a tensor in the form of tensor objects without referencing them to any basis. In the basis representation the derivatives of the tensor argument functions are given in works [15-17].

The rules for constructing derivatives of the scalar and tensor functions of the tensor argument with respect to a tensor and the form of their representation depend on the accepted scheme of the interaction of the basis vectors in a double scalar product of two tensors. For the second-rank tensors, there are two such independent schemes, for the second and fourth-rank tensors there are 12 schemes and for tensors of the fourth rank their number increases to 36. As a result, there are no uniform rules and uniform form of writing derivatives of the scalar and tensor functions of the tensor argument with respect to a tensor. Until the present time there have been only three groups of the rules and forms of derivative representation. Their derivation and specific expressions for derivatives for each group are presented in works [8-10], [11-13] and [15-18], respectively. It is interesting to note that in contrast to [15-17] of the last group of articles, in which the derivatives of the tensor argument functions are written in the basic form, work [18] uses the concept of positional scalar multiplication [19], which allowed the authors to represent the rules for constructing derivatives and the derivatives themselves in the form of tensor objects without referencing them to any basis. In work [14], the correspondence between the results obtained in the works of the first and second group has been established. However, the results obtained in the works of the third group, for some reason, were not taken into account in this consideration.

The rules for constructing derivatives of the scalar and tensor functions of the tensor argument with respect to a tensor in the works of the first group is based on a widely used scheme of the interaction of the basis vectors in a double scalar product of the second and fourth-rank tensors (the schemes of the basis vector interaction will be considered in the next Section). However, this yields rather exotic representation forms, which essentially differ from the forms of derivative representation commonly used in conventional mathematical analysis.

In the works of the second group the rules of constructing derivatives for the scalar and tensor functions of the tensor argument with respect to a tensor and the forms of derivative representation do not differ from the rules of derivative construction and the forms of their representation used in conventional mathematical

*For a second-rank tensor, the transposition operation reduces to a permutation of the basis vectors of this tensor and is denoted by the superscript *T*.

analysis. However, it must be emphasized that their construction is based on a special (given below) scheme of the basis vector interaction in a double scalar product of the second and fourth-rank tensors.

Finally, the work of the third group also operates with a widespread scheme of the interaction of the basis vectors, which however is different from those used in the works of the first group. Furthermore, the resulting rules for constructing derivatives of the scalar and tensor functions of the tensor argument with respect to a tensor and the forms of derivative representation differ insignificantly from the corresponding rules and forms of representation used in conventional mathematical analysis.

Note that derivatives of the scalar functions of the tensor argument with respect to a tensor are similar in any group. The only difference is in the derivatives of the tensor functions of the tensor arguments. Bearing this in mind, we may suppose that such ambiguity will not be critical to the analysis of solid mechanics who are one of the main users of the results of this article.

## II. Double Scalar Product of Tensors

The rules for constructing derivatives of the scalar and tensor functions of the tensor argument with respect to a tensor and the forms of representation of the derivatives themselves depend on the scheme of double scalar multiplication of the basis vectors of the second and fourth-rank tensors, which lies at the heart of the theory that allows you to formulate these rules. For the second-rank tensors $\mathbf{A}$ and $\mathbf{B}$, there are only two independent schemes of double scalar multiplication of the basis vectors, which can be represented as follows:

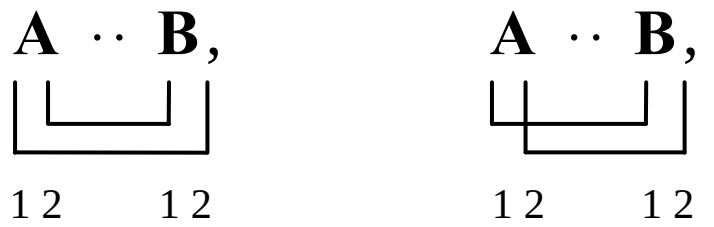

(2.1)

Here, the horizontal brackets indicate for which basis vectors and of which tensors the scalar multiplication is carried out. For example, in the left scheme the second basis vector of the second-rank tensor $\mathbf{A}$ is multiplied by the first basis vector of the second-rank tensor $\mathbf{B}$, while the first basis vector of the tensor $\mathbf{A}$ is multiplied by the second basis vector of the tensor $\mathbf{B}$ (first, the nearest vectors are multiplied and turn into a scalar, then, multiplication is performed for the next vectors, which have now become the nearest). The right-hand multiplication scheme is also of frequent use. It is often referred to as a double scalar product and its identification is also based on two points (horizontal or vertical). The cross multiplication of the nearest basis vectors is fulfilled according to this scheme. We will denote this operation by the symbol $(\cdot\cdot)$, $\mathbf{A}(\cdot\cdot)\mathbf{B}$, leaving the operation notation by two horizontal (or vertical) points $\cdot\cdot$ for the left scheme of the basis vector multiplication. Certainly,

$$\mathbf{A}(\cdot\cdot)\mathbf{B} = \mathbf{A}\cdot\cdot\mathbf{B}^T = \mathbf{A}^T\cdot\cdot\mathbf{B}, \tag{2.2}$$

i.e., the operations of cross-multiplication of the nearest basic vectors and their successive multiplication are interrelated.

Schemes that are similar to the above mentioned ones and exhibiting the relationship between the basis vectors at the interaction of tensors are of demonstrative character. They allow us to represent the operations with tensor objects without referencing the latter to any basis and will be used later in the article. In particular, based on the schemes (2.1), we can easily demonstrate the validity of the following relations:

$$\mathbf{A} * \mathbf{B} = \mathbf{B} * \mathbf{A} = \mathbf{A}^T * \mathbf{B}^T = \mathbf{B}^T * \mathbf{A}^T, \tag{2.3}$$

where $* = \cdot\cdot$ or $* = (\cdot\cdot)$. Also, from (2.1) – (2.3), we obtain the following relationships:

$$\mathbf{A}\cdot\cdot(\mathbf{B}\cdot\mathbf{C}) = (\mathbf{A}\cdot\mathbf{B})\cdot\cdot\mathbf{C}, \qquad \mathbf{A}(\cdot\cdot)(\mathbf{B}\cdot\mathbf{C}) = (\mathbf{A}^T\cdot\mathbf{B})(\cdot\cdot)\mathbf{C}^T, \tag{2.4}$$

the first of which is more attractive because of its simplicity and clarity (here we need to permute only the operation $\cdot\cdot$, in contrast to the second relationship, where in addition to permutation of the $(\cdot\cdot)$ operation it is necessary to perform transposition of certain tensors).

As noted above, for the second-rank tensors there are only two independent schemes of double scalar multiplication of the basis vectors, which are shown in (2.1). For each of the schemes we introduced the appropriate operation sign. The double scalar multiplication of the basis vectors of the second and fourth-rank tensors involves already 12 independent schemes of vector interaction, so that there is no use of introducing the appropriate operation sign for each of them. Therefore, to identify the operation of the scalar multiplication we use in the general case a point with an indication of the basis vectors of the left (above the point) and the right (beneath the point) tensors of arbitrary ranks interrelated by this operation: $\overset{m}{\underset{n}{\cdot}}$. Then, in the accepted notation, the relations (2.1) are given as $\mathbf{A}\cdot\cdot\mathbf{B} = \mathbf{A} \mathop{\cdot\cdot}\limits^{1\,2}_{2\,1} \mathbf{B}$ and $\mathbf{A}(\cdot\cdot)\mathbf{B} = \mathbf{A} \mathop{\cdot\cdot}\limits^{1\,2}_{1\,2} \mathbf{B}$ from which we immediately arrive at the expressions (2.2).

Of the twelve independent schemes of the interaction between the basis vectors of the second-rank $\mathbf{A}$ and fourth-rank $\mathbf{H}^{IV}$ tensors during the operation of their double multiplication, the schemes $\mathbf{A} \mathop{\cdot\cdot}\limits^{1\,2}_{1\,4} \mathbf{H}^{IV}$ and $\mathbf{H}^{IV} \mathop{\cdot\cdot}\limits^{2\,3}_{1\,2} \mathbf{A}$ are introduced and used in the literature for constructing the derivatives of the scalar and tensor functions of the tensor argument with respect to this argument. Below, these schemes are represented in detail as

$$\mathbf{A}\cdot\cdot\,\mathbf{H}^{IV}, \qquad \mathbf{H}^{IV}\,\cdot\cdot\,\mathbf{A}, \tag{2.5}$$

12 1234 1 2 3 4 1 2

as well as the schemes $\mathbf{A} \mathop{\cdot\cdot}\limits^{1\,2}_{1\,2} \mathbf{H}^{IV}$ and $\mathbf{H}^{IV} \mathop{\cdot\cdot}\limits^{3\,4}_{1\,2} \mathbf{A},$ which are represented in the following visual form:

$$\mathbf{A}\cdot\cdot\,\mathbf{H}^{IV}, \qquad \mathbf{H}^{IV}\,\cdot\cdot\,\mathbf{A}, \tag{2.6}$$

12 1 2 34 1 2 3 4 1 2

In the last line of these schemes, the basis vectors of each tensor are numbered from the left to the right. Unlike the schemes (2.5) and (2.6), the operation of obtaining the double scalar product introduced above for the second-rank tensors, and implying a sequential scalar multiplication of the nearest basis vectors, will also be used for the second- and fourth-rank tensors: $\mathbf{A} \mathop{\cdot\cdot}\limits^{1\,2}_{2\,1} \mathbf{H}^{IV}$ and $\mathbf{H}^{IV} \mathop{\cdot\cdot}\limits^{3\,4}_{2\,1} \mathbf{A}.$ For the sake of clarity these schemes are presented below at greater length:

$$\mathbf{A}\cdot\cdot\,\mathbf{H}^{IV}, \qquad \mathbf{H}^{IV}\,\cdot\cdot\,\mathbf{A}, \tag{2.7}$$

12 12 34 1 2 3 4 1 2

For a double multiplication of the basis vectors of the two fourth-rank tensors $\mathbf{P}^{IV}$ and $\mathbf{H}^{IV},$ there are already 36 independent schemes of the interaction. The number of independent schemes of the basis vector interaction in the double multiplication of the basis vectors $\gamma$ of the two $M-$ and $N-$rank tensors is defined by the ratio $\gamma = a(b-1)(a-1),$ where $a = \max(M,N)$ and $b = \min(M,N),$ and $\gamma = 2$ when $M = N = 2,$ $\gamma = 12$ at $M = 2, N = 4$ and $\gamma = 36$ at $M = N = 4,$ as it has been noted above. Of these 36 schemes of interaction, only three are used for constructing the derivatives of the scalar and tensor functions of the tensor argument with respect to a tensor: $\mathbf{P}^{IV} \mathop{\cdot\cdot}\limits^{2\,3}_{1\,4} \mathbf{H}^{IV},$ $\mathbf{P}^{IV} \mathop{\cdot\cdot}\limits^{3\,4}_{1\,2} \mathbf{H}^{IV},$ $\mathbf{P}^{IV} \mathop{\cdot\cdot}\limits^{3\,4}_{2\,1} \mathbf{H}^{IV},$

$$\mathbf{P}^{IV}\cdot\cdot\,\mathbf{H}^{IV}, \qquad \mathbf{P}^{IV}\cdot\cdot\,\mathbf{H}^{IV}, \qquad \mathbf{P}^{IV}\cdot\cdot\,\mathbf{H}^{IV}. \tag{2.8}$$

1234 1234 1234 1234 1234 1234

Schemes (2.5) and the first scheme in (2.8) were used in [11-13], schemes (2.6) and the second one in (2.8) are used in [8-10], while scheme (2.7) and the last scheme in (2.8) – in [15-18] and in this article.

In scheme (2.6) and in the second scheme (2.8) the operation of cross multiplication of the nearest basis vectors of the tensors is used. For this operation, we have introduced the notation $(\cdot\cdot)$ and used it to rewrite the above schemes as follows:

$$\mathbf{A}(\cdot\cdot)\mathbf{H}^{IV}, \qquad \mathbf{H}^{IV}(\cdot\cdot)\mathbf{A}, \qquad \mathbf{P}^{IV}(\cdot\cdot)\mathbf{H}^{IV}. \tag{2.9}$$

In scheme (2.7) and the last scheme (2.8) the sequential multiplication of the nearest basis vectors of the tensors is employed. For this operation we have introduced the notation $\cdot\cdot$ and used it to rewrite the above schemes as follows:

$$\mathbf{A}\cdot\cdot\mathbf{H}^{IV}, \qquad \mathbf{H}^{IV}\cdot\cdot\mathbf{A}, \qquad \mathbf{P}^{IV}\cdot\cdot\mathbf{H}^{IV}. \tag{2.10}$$

Finally, after introducing the operation $[\cdot\cdot]$, which is called the operation of double positional scalar multiplication and involves multiplication of the basis vectors at the interaction of the second and fourth-rank tensors, the fourth- and second-rank tensors and fourth-rank tensors in accordance with scheme (2.5) and the first scheme in (2.8), we can represent these schemes as

$$\mathbf{A}[\cdot\cdot]\mathbf{H}^{IV}, \qquad \mathbf{H}^{IV}[\cdot\cdot]\mathbf{A}, \qquad \mathbf{P}^{IV}[\cdot\cdot]\mathbf{H}^{IV}. \tag{2.11}$$

The substance of this operation can easily be interpreted as follows: the tenor $\mathbf{B}$ of the second or fourth rank on the right-hand side of relations is substituted in the form of $(\cdot\mathbf{B}\cdot)$ in the middle of the group of the basis vectors of the left tensor $\mathbf{A}$, which is also of the second or fourth rank, and the operations of scalar multiplication specified for tensor $\mathbf{B}$ are performed. For example, for the second-rank tensors $\mathbf{A}$ and $\mathbf{B}$, which are represented, respectively by contravariant and covariant components in the basis $\mathbf{r}_i$ and $\mathbf{r}^i$ (see the beginning of the next section), we have

$$\mathbf{A}[\cdot\cdot]\mathbf{B} = A^{ij}\mathbf{r}_i\mathbf{r}_j[\cdot\cdot]B_{kp}\mathbf{r}^k\mathbf{r}^p = A^{ij}(\mathbf{r}_i\cdot\mathbf{B}\cdot\mathbf{r}_j) = A^{ij}(\mathbf{r}_i\cdot B_{kp}\mathbf{r}^k\mathbf{r}^p\cdot\mathbf{r}_j) = A^{ij}B_{ij} = \mathbf{A}(\cdot\cdot)\mathbf{B}, \tag{2.12}$$

that is, for tensors of the second rank the operations $[\cdot\cdot]$ and $(\cdot\cdot)$ do not differ. For the second-rank tensor $\mathbf{A}$ and the fourth-rank tensor $\mathbf{H}^{IV}$ we obtain

$$\mathbf{A}[\cdot\cdot]\mathbf{H}^{IV} = A^{ij}\mathbf{r}_i\mathbf{r}_j[\cdot\cdot]H_{mnkp}\mathbf{r}^m\mathbf{r}^n\mathbf{r}^k\mathbf{r}^p = A^{ij}(\mathbf{r}_i\cdot\mathbf{H}^{IV}\cdot\mathbf{r}_j) = A^{ij}(\mathbf{r}_i\cdot H_{mnkp}\mathbf{r}^m\mathbf{r}^n\mathbf{r}^k\mathbf{r}^p\cdot\mathbf{r}_j) = A^{ij}H_{inkj}\mathbf{r}^n\mathbf{r}^k,$$

$$\mathbf{H}^{IV}[\cdot\cdot]\mathbf{A} = H_{ijkp}\mathbf{r}^i\mathbf{r}^j\mathbf{r}^k\mathbf{r}^p[\cdot\cdot]A^{mn}\mathbf{r}_m\mathbf{r}_n = H_{ijkp}(\mathbf{r}^i\mathbf{r}^j\cdot\mathbf{A}\cdot\mathbf{r}^k\mathbf{r}^p) =$$
$$= H_{ijkp}(\mathbf{r}^i\mathbf{r}^j\cdot A^{mn}\mathbf{r}_m\mathbf{r}_n\cdot\mathbf{r}^k\mathbf{r}^p) = H_{ijkp}A^{jk}\mathbf{r}^i\mathbf{r}^p,$$

which corresponds to the schemes (2.5). Finally, for the fourth-rank tensors $\mathbf{P}^{IV}$ and $\mathbf{H}^{IV}$ we have

$$\mathbf{P}^{IV}[\cdot\cdot]\mathbf{H}^{IV} = P^{ijkl}\mathbf{r}_i\mathbf{r}_j\mathbf{r}_k\mathbf{r}_l[\cdot\cdot]H_{mnsp}\mathbf{r}^m\mathbf{r}^n\mathbf{r}^s\mathbf{r}^p = P^{ijkl}(\mathbf{r}_i\mathbf{r}_j\cdot\mathbf{H}^{IV}\cdot\mathbf{r}_k\mathbf{r}_l) =$$
$$= P^{ijkl}(\mathbf{r}_i\mathbf{r}_j\cdot H_{mnsp}\mathbf{r}^m\mathbf{r}^n\mathbf{r}^s\mathbf{r}^p\cdot\mathbf{r}_k\mathbf{r}_l) = P^{ijkl}H_{jnsk}\mathbf{r}_i\mathbf{r}^n\mathbf{r}^s\mathbf{r}_l,$$

which corresponds to the first scheme (2.8) and completes its definition by specifying the sequence order of the noninteracting basis vectors in the resulting tensor. Note that in [14] the operation $[\cdot\cdot]$ is designated as $\bullet\circ$.

In what follows, we will use the above introduced three groups of relations between tensors (2.9) – (2.11).

## III. Fourth-Rank Isotropic Tensors and Their Properties

In view of the fact that the solid mechanics is one of the main fields of research, for which the results of this article might be beneficial, we will build the vector and tensor objects in the basis sets, which are directly related to the initial $\kappa_0$, current $\kappa$ and any intermediate $\kappa_*$ configurations, which the body can occupy in the process of deformation. The position of a material point in the configuration $\kappa_0$ is determined by the radius vector $\mathbf{r}(q^i)$, in the configuration $\kappa$ – by the radius vector $\mathbf{R}(q^i)$ and in the configuration $\kappa_*$ – by the radius vector $\breve{\mathbf{R}}(q^i)$, where $q^i$, $i = 1, 2, 3$ are the generalized Lagrangian (material) coordinates. Then, the vectors of the main basis in these configurations are

$$\mathbf{r}_i = \partial\mathbf{r}/\partial q^i, \qquad \mathbf{R}_i = \partial\mathbf{R}/\partial q^i, \qquad \breve{\mathbf{R}}_i = \partial\breve{\mathbf{R}}/\partial q^i,$$

the reciprocal basis vectors $\mathbf{r}^i$, $\mathbf{R}^i$, $\breve{\mathbf{R}}^i$ are found from the conditions

$$\mathbf{r}_i\cdot\mathbf{r}^j = \delta_i^j, \qquad \mathbf{R}_i\cdot\mathbf{R}^j = \delta_i^j, \qquad \breve{\mathbf{R}}_i\cdot\breve{\mathbf{R}}^j = \delta_i^j,$$

where $\delta_i^j$ is the Kronecker Delta and the simplest expressions for the metric tensors of these configurations $\mathbf{g}$, $\mathbf{G}$, $\breve{\mathbf{G}}$ can be written as

$$\mathbf{g} = \mathbf{r}_i\mathbf{r}^i = \mathbf{r}^i\mathbf{r}_i, \qquad \mathbf{G} = \mathbf{R}_i\mathbf{R}^i = \mathbf{R}^i\mathbf{R}_i, \qquad \breve{\mathbf{G}} = \breve{\mathbf{R}}_i\breve{\mathbf{R}}^i = \breve{\mathbf{R}}^i\breve{\mathbf{R}}_i. \tag{3.1}$$

Here, for a second-rank tensor, instead of the tensor product operation $\mathbf{a} \otimes \mathbf{b}$ for two vectors, we used a more laconic form of writing, namely, the dyad representation $\mathbf{ab}$, and applied the convention of summation over the repeated indices (Einstein summation convention). In the following, we will use this form of representation.

We will use three fourth-rank isotropic tensors $\mathbf{C}_I^{IV}$, $\mathbf{C}_{II}^{IV}$ and $\mathbf{C}_{III}^{IV}$. The first of them is defined as an ordinary tensor product of two metric tensors of the initial $\mathbf{g}$, current $\mathbf{G}$ or any intermediate $\breve{\mathbf{G}}$ configuration. In their configurations these metric tensors play the role of the unit tensors and, by virtue of their equality as the tensor objects, will be denoted by $\mathbf{I}$ if there is no need to specify the basis of the configuration, to which the tensors refer. In particular, the first isotropic tensor of the fourth rank is represented in the basis of the initial configuration as

$$\mathbf{C}_I^{IV} = \mathbf{gg} = \mathbf{g} \otimes \mathbf{g}$$

or, with account of (1.1) and (3.1), as

$$\mathbf{C}_I^{IV} = \mathbf{g} \otimes \mathbf{g} = \mathbf{r}_i \mathbf{r}^i \mathbf{r}_j \mathbf{r}^j = \mathbf{r}_i \mathbf{r}^i \mathbf{r}^j \mathbf{r}_j = \mathbf{r}^i \mathbf{r}_i \mathbf{r}_j \mathbf{r}^j = \mathbf{r}^i \mathbf{r}_i \mathbf{r}^j \mathbf{r}_j, \tag{3.2}$$

where for each tensor $\mathbf{g}$ the first or second representation from (3.1) is used.

Using the operations $\boxtimes$ and $\hat{\boxtimes}$ and any representation of the tensor $\mathbf{g}$ from (3.1), we introduce another two isotropic tensors of the fourth rank:

$$\mathbf{C}_{II}^{IV} = \mathbf{g} \boxtimes \mathbf{g} = \mathbf{r}_i \mathbf{r}_j \mathbf{r}^i \mathbf{r}^j = \mathbf{r}_i \mathbf{r}^j \mathbf{r}^i \mathbf{r}_j = \mathbf{r}^i \mathbf{r}_j \mathbf{r}_i \mathbf{r}^j = \mathbf{r}^i \mathbf{r}^j \mathbf{r}_i \mathbf{r}_j \tag{3.3}$$

and

$$\mathbf{C}_{III}^{IV} = \mathbf{g} \hat{\boxtimes} \mathbf{g} = \mathbf{r}_i \mathbf{r}_j \mathbf{r}^j \mathbf{r}^i = \mathbf{r}_i \mathbf{r}^j \mathbf{r}_j \mathbf{r}^i = \mathbf{r}^i \mathbf{r}_j \mathbf{r}^j \mathbf{r}_i = \mathbf{r}^i \mathbf{r}^j \mathbf{r}_j \mathbf{r}_i. \tag{3.4}$$

These fourth-rank tensors are called isotropic because, as noted above, they are constructed in the basis of tensors $\mathbf{g}$, $\mathbf{G}$ or $\breve{\mathbf{G}}$, which, if represented in the orthonormal (Cartesian) basis $\mathbf{i}_i$, $i = 1, 2, 3$, have the same form of representation (coincide)

$$\mathbf{g} = \mathbf{G} = \breve{\mathbf{G}} = \mathbf{I} = \mathbf{i}_i \mathbf{i}_i \tag{3.5}$$

and are isotropic, since their parity group (equvalence group) is a complete orthogonal group: $\mathbf{O} \cdot \mathbf{I} \cdot \mathbf{O}^T = \mathbf{I}$ for any orthogonal tensor $\mathbf{O}$.

The properties of the isotropic fourth-rank tensors are derived from relations (3.2) – (3.4) and depend on which of the schemes of interaction between the basis vectors of the second-rank tensor $\mathbf{A}$ and fourth-rank tensor $\mathbf{H}^{IV}$ is used in the operation of their double multiplication. These properties are given in Table 1.

**Table 1** *The properties of the isotropic fourth-rank tensors*

| Scheme | Operation | $\mathbf{C}_I^{IV}$ | $\mathbf{C}_{II}^{IV}$ | $\mathbf{C}_{III}^{IV}$ | References |
|---|---|---|---|---|---|
| (2.9) | $\mathbf{A}(\cdot\cdot)\mathbf{H}^{IV}$, $\mathbf{H}^{IV}(\cdot\cdot)\mathbf{A}$ | $I_1(\mathbf{A})\mathbf{g}$ | $\mathbf{A} \Rightarrow \mathbf{C}_{II}^{IV} \equiv \mathbb{I}$ | $\mathbf{A}^T \Rightarrow \mathbf{C}_{III}^{IV} \equiv \mathbb{T}$ | [8-10] |
| (2.10) | $\mathbf{A} \cdot\cdot \mathbf{H}^{IV}$, $\mathbf{H}^{IV} \cdot\cdot \mathbf{A}$ | $I_1(\mathbf{A})\mathbf{g}$ | $\mathbf{A}^T \Rightarrow \mathbf{C}_{II}^{IV} \equiv \mathbb{T}$ | $\mathbf{A} \Rightarrow \mathbf{C}_{III}^{IV} \equiv \mathbb{I}$ | [15-18] |
| (2.11) | $\mathbf{A}[\cdot\cdot]\mathbf{H}^{IV}$, $\mathbf{H}^{IV}[\cdot\cdot]\mathbf{A}$ | $\mathbf{A} \Rightarrow \mathbf{C}_I^{IV} \equiv \mathbb{I}$ | $\mathbf{A}^T \Rightarrow \mathbf{C}_{II}^{IV} \equiv \mathbb{T}$ | $I_1(\mathbf{A})\mathbf{g}$ | [11-13] |

The cell of the table at the intersection of the row and column comprises the result of the operation specified in the second column of the row, in which the tensor $\mathbf{H}^{IV}$ is replaced by the corresponding isotropic tensor of the fourth rank (3.2) – (3.4) specified in columns 3 – 5. Scalar multiplication of tensors presented in Table 1 can easily be realized by taking into account the orthogonality of the main and reciprocal basis vectors. To this end, it is necessary to choose from the four forms of representation of each of the fourth-rank tensors the form that will correspond to the form of representation of the second-rank tensor $\mathbf{A}$ (in the bases $\mathbf{r}_i$ and $\mathbf{r}^i$ this tensor has four forms of representation – contravariant, covariant and two mixed forms). For example, if tensor $\mathbf{A}$ is written as $\mathbf{A} = A^{ij}\mathbf{r}_i\mathbf{r}_j$, then in the products $\mathbf{A}(\cdot\cdot)\mathbf{C}_{II}^{IV}$ and $\mathbf{A} \cdot\cdot \mathbf{C}_{II}^{IV}$ (the second and third rows of the Table 1) the tensor $\mathbf{C}_{II}^{IV}$ is convenient to use in the form of its last representation (3.3), whereas in the products $\mathbf{C}_{II}^{IV}(\cdot\cdot)\mathbf{A}$ and $\mathbf{C}_{II}^{IV} \cdot\cdot \mathbf{A}$ the first form of representation is more suitable.

Table 1 also presents the consequences arising from the obtained results – the equivalence of the tensors $\mathbb{I}$ and $\mathbb{T}$ to the fourth-rank tensors (3.2) – (3.4) in the double scalar multiplication of the latter by a second rank tensor realized according to a certain scheme of the basic vector interaction. The properties of the tensors $\mathbb{I}$ and $\mathbb{T}$ are described by relations (1.2). From the second and third rows of Table 1 the relation $(\cdot\cdot) = \cdot\cdot\,\mathbf{C}_{II}^{IV}\,\cdot\cdot$ follows, from which we readily obtain the expressions (2.2) – (2.4) for the second-rank tensors. Indeed,

$$\mathbf{A}(\cdot\cdot)\mathbf{B} = \mathbf{A}\cdot\cdot\mathbf{C}_{II}^{IV}\cdot\cdot\mathbf{B} = \begin{cases} (\mathbf{A}\cdot\cdot\mathbf{C}_{II}^{IV})\cdot\cdot\mathbf{B} = \mathbf{A}^T\cdot\cdot\mathbf{B} \\ \mathbf{A}\cdot\cdot(\mathbf{C}_{II}^{IV}\cdot\cdot\mathbf{B}) = \mathbf{A}\cdot\cdot\mathbf{B}^T \end{cases} \quad \Rightarrow \qquad (2.2)_1$$

$$\Rightarrow \qquad \mathbf{A}(\cdot\cdot)\mathbf{B} = \begin{cases} \mathbf{A}^T\cdot\cdot\mathbf{B} = \mathbf{A}^T\cdot\cdot(\mathbf{C}_{II}^{IV}\cdot\cdot\mathbf{B}^T) = \mathbf{A}^T(\cdot\cdot)\mathbf{B}^T \\ \mathbf{B}\cdot\cdot\mathbf{A}^T = \begin{cases} \mathbf{B}\cdot\cdot(\mathbf{C}_{II}^{IV}\cdot\cdot\mathbf{A}) = \mathbf{B}(\cdot\cdot)\mathbf{A} \\ (\mathbf{B}^T\cdot\cdot\mathbf{C}_{II}^{IV})\cdot\cdot\mathbf{A}^T = \mathbf{B}^T(\cdot\cdot)\mathbf{A}^T \end{cases} \end{cases} ; \qquad (2.3)_1$$

$$\mathbf{A}(\cdot\cdot)(\mathbf{B}\cdot\mathbf{C}) = \mathbf{A}\cdot\cdot\mathbf{C}_{II}^{IV}\cdot\cdot(\mathbf{B}\cdot\mathbf{C}) = (\mathbf{A}\cdot\cdot\mathbf{C}_{II}^{IV})\cdot\cdot(\mathbf{B}\cdot\mathbf{C}) = \mathbf{A}^T\cdot\cdot(\mathbf{B}\cdot\mathbf{C}) =$$
$$= (\mathbf{A}^T\cdot\mathbf{B})\cdot\cdot\mathbf{C} = (\mathbf{A}^T\cdot\mathbf{B})\cdot\cdot(\mathbf{C}_{II}^{IV}\cdot\cdot\mathbf{C}^T) = (\mathbf{A}^T\cdot\mathbf{B})(\cdot\cdot)\mathbf{C}^T. \qquad (2.4)_1$$

Using the relation $(\cdot\cdot) = \cdot\cdot\,\mathbf{C}_{II}^{IV}\,\cdot\cdot$, we obtain from relations (2.9) and (2.10) the following expressions for the second- and fourth-rank tensors:

$$\mathbf{A}(\cdot\cdot)\mathbf{H}^{IV} = \mathbf{A}\cdot\cdot\mathbf{C}_{II}^{IV}\cdot\cdot\mathbf{H}^{IV} = (\mathbf{A}\cdot\cdot\mathbf{C}_{II}^{IV})\cdot\cdot\mathbf{H}^{IV} = \mathbf{A}^T\cdot\cdot\mathbf{H}^{IV},$$
$$\mathbf{H}^{IV}(\cdot\cdot)\mathbf{A} = \mathbf{H}^{IV}\cdot\cdot\mathbf{C}_{II}^{IV}\cdot\cdot\mathbf{A} = \mathbf{H}^{IV}\cdot\cdot(\mathbf{C}_{II}^{IV}\cdot\cdot\mathbf{A}) = \mathbf{H}^{IV}\cdot\cdot\mathbf{A}^T,$$

while for tensors of the fourth rank we have

$$\mathbf{P}^{IV}(\cdot\cdot)\,\mathbf{H}^{IV} = \mathbf{P}^{IV}\cdot\cdot\mathbf{C}_{II}^{IV}\cdot\cdot\,\mathbf{H}^{IV} = \begin{cases} (\mathbf{P}^{IV}\cdot\cdot\mathbf{C}_{II}^{IV})\cdot\cdot\,\mathbf{H}^{IV} = (\mathbf{P}^{IV})^{dr}\cdot\cdot\,\mathbf{H}^{IV} \\ \mathbf{P}^{IV}\cdot\cdot(\mathbf{C}_{II}^{IV}\cdot\cdot\,\mathbf{H}^{IV}) = \mathbf{P}^{IV}\cdot\cdot(\mathbf{H}^{IV})^{dl} \end{cases},$$

where $(\mathbf{M}^{IV})^{dr}$ and $(\mathbf{M}^{IV})^{dl}$ are two of the variants of transpose of the fourth-rank tensor $\mathbf{M}^{IV}$ introduced in [14]. These variants correspond to the operations $(\mathbf{M}^{IV})^{dr} = \mathbf{M}^{IV}\cdot\cdot\mathbf{C}_{II}^{IV}$ and $(\mathbf{M}^{IV})^{dl} = \mathbf{C}_{II}^{IV}\cdot\cdot\,\mathbf{M}^{IV}$ leading to the following permutation of the basis vectors of the tensor $\mathbf{M}^{IV}$:

$$\underset{1234}{(\mathbf{M}^{IV})^{dr}} = \underset{1243}{\mathbf{M}^{IV}}, \qquad \underset{1234}{(\mathbf{M}^{IV})^{dl}} = \underset{2134}{\mathbf{M}^{IV}}.$$

In conclusion, we note that in this Section the metric tensors of the initial configuration have been used to construct the isotropic tensors of the fourth rank in relations (3.2) – (3.4). As mentioned above, these tensors coincide as the tensor objects with tensors $\mathbf{G}, \breve{\mathbf{G}}$ and $\mathbf{I}$: $\mathbf{g} = \mathbf{G} = \breve{\mathbf{G}} = \mathbf{I}$ (see (3.5)). Therefore, in the current configuration, the isotropic fourth-rank tensors $\mathbf{C}_I^{IV}$, $\mathbf{C}_{II}^{IV}$ and $\mathbf{C}_{III}^{IV}$ are represented as

$$\mathbf{C}_I^{IV} = \mathbf{R}_i\mathbf{R}^i\mathbf{R}_j\mathbf{R}^j = \mathbf{R}_i\mathbf{R}^i\mathbf{R}^j\mathbf{R}_j = \mathbf{R}^i\mathbf{R}_i\mathbf{R}_j\mathbf{R}^j = \mathbf{R}^i\mathbf{R}_i\mathbf{R}^j\mathbf{R}_j,$$
$$\mathbf{C}_{II}^{IV} = \mathbf{R}_i\mathbf{R}_j\mathbf{R}^i\mathbf{R}^j = \mathbf{R}_i\mathbf{R}^j\mathbf{R}^i\mathbf{R}_j = \mathbf{R}^i\mathbf{R}_j\mathbf{R}_i\mathbf{R}^j = \mathbf{R}^i\mathbf{R}^j\mathbf{R}_i\mathbf{R}_j$$

and

$$\mathbf{C}_{III}^{IV} = \mathbf{R}_i\mathbf{R}_j\mathbf{R}^j\mathbf{R}^i = \mathbf{R}_i\mathbf{R}^j\mathbf{R}_j\mathbf{R}^i = \mathbf{R}^i\mathbf{R}_j\mathbf{R}^j\mathbf{R}_i = \mathbf{R}^i\mathbf{R}^j\mathbf{R}_j\mathbf{R}_i$$

and in the orthonormal Cartesian basis, in which the main and reciprocal reference vectors coincide, they take the following form:

$$\mathbf{C}_I^{IV} = \mathbf{i}_i\mathbf{i}_i\mathbf{i}_j\mathbf{i}_j, \qquad \mathbf{C}_{II}^{IV} = \mathbf{i}_i\mathbf{i}_j\mathbf{i}_i\mathbf{i}_j, \qquad \mathbf{C}_{III}^{IV} = \mathbf{i}_i\mathbf{i}_j\mathbf{i}_j\mathbf{i}_i. \qquad (3.6)$$

## IV. Differentiation of Scalar Functions of a Tensor Argument

Let the scalar function $\Phi$ depend on the vector argument $\mathbf{a}$, $\Phi = \Phi(\mathbf{a})$, which can be represented by its contravariant and covariant components:

$$\mathbf{a} = a^i\,\mathbf{r}_i = a_i\,\mathbf{r}^i.$$

Here $\mathbf{r}_i$ is the main, generally not orthonormal basis, and $\mathbf{r}^i$ is its reciprocal. A change (variation, increment) of the scalar function $\Phi$ is determined by a change (variation) of the vector argument $\mathbf{a}$, which occurs both due to the variation of the coordinate components, and, in the general case, due to the variation of the basis vectors:

$d\mathbf{a}=(da^i)\mathbf{r}_i+a^i\,d\mathbf{r}_i=(da_i)\mathbf{r}^i+a_i\,d\mathbf{r}^i$. The representation of the vector $\mathbf{a}$ in the Cartesian basis $\mathbf{i}_i$, which does not change neither in the magnitude nor in the direction, $\mathbf{a}=\hat{a}^i\,\mathbf{i}_i$, makes it possible to determine a change of the vector $\mathbf{a}$ only by changes in its coordinate components $\hat{a}^i$: $d\mathbf{a}=(d\hat{a}^i)\,\mathbf{i}_i$. So $\Phi(\mathbf{a})$ is the function of $\hat{a}^i$, $\Phi(\mathbf{a})=\Phi(\hat{a}^i)$, only. Then, the conventional mathematical analysis of the function of few variables (in this case, three) yields

$$d\,\Phi(\mathbf{a})=\frac{\partial\Phi(\hat{a}^1,\,\hat{a}^2,\,\hat{a}^3)}{\partial\,\hat{a}^i}\,d\,\hat{a}^i,\quad i=1,2,3. \tag{4.1}$$

Using a particular case of the inverse tensor attribute theorem* (see [20]), this expression can be rewritten as

$$d\,\Phi(\mathbf{a})=\frac{\partial\Phi(\hat{a}^1,\,\hat{a}^2,\,\hat{a}^3)}{\partial\,\hat{a}^i}\,\mathbf{i}_i\cdot(d\,\hat{a}^j)\,\mathbf{i}_j \tag{4.2}$$

or, given that $(d\hat{a}^i)\,\mathbf{i}_i=d\mathbf{a}$ as

$$d\,\Phi(\mathbf{a})=\frac{\partial\Phi(\hat{a}^1,\,\hat{a}^2,\,\hat{a}^3)}{\partial\,\hat{a}^i}\,\mathbf{i}_i\cdot d\,\mathbf{a}=\frac{\partial\,\Phi(\mathbf{a})}{\partial\,\mathbf{a}}\cdot d\,\mathbf{a}, \tag{4.3}$$

where

$$\frac{\partial\,\Phi(\mathbf{a})}{\partial\,\mathbf{a}}=\frac{\partial\,\Phi(\hat{a}^1,\,\hat{a}^2,\,\hat{a}^3)}{\partial\,\hat{a}^i}\,\mathbf{i}_i. \tag{4.4}$$

Knowing the derivative $\partial\,\Phi(\mathbf{a})/\partial\,\mathbf{a}$ in the Cartesian basis, we can easily represent it in any other system of coordinates.

Now let the scalar function $\Phi$ depend on the tensor of the second rank $\mathbf{A}$, $\Phi=\Phi(\mathbf{A})$, which can be represented both in the variable basis sets $\mathbf{r}_i$, $\mathbf{r}^i$, and in the invariable Cartesian basis $\mathbf{i}_i:\mathbf{A}=a^{ij}\mathbf{i}_i\mathbf{i}_j$. In the latter case $d\mathbf{A}=(da^{ij})\,\mathbf{i}_i\mathbf{i}_j$, $\Phi(\mathbf{A})=\Phi(a^{11},a^{12},...,a^{33})$ and a generalization of expression (4.1) is written as

$$d\Phi(\mathbf{A})=\frac{\partial\Phi(a^{11},\,a^{12},\ldots,\,a^{33})}{\partial a^{i\,j}}\;d\,a^{i\,j},\quad i,j=1,2,3, \tag{4.5}$$

which, using a special case of the inverse tensor attribute theorem [20] (see footnote on this page), can be represented in two forms according to the schemes (2.1):

$$d\,\Phi(\mathbf{A})=\frac{\partial\Phi(a^{11},\,a^{12},\ldots,\,a^{33})}{\partial\,a^{i\,j}}\,\mathbf{i}_i\,\mathbf{i}_j\cdot\cdot\,(d\,a^{k\,l})\,\mathbf{i}_l\;\mathbf{i}_k$$

or

$$d\,\Phi(\mathbf{A})=\frac{\partial\Phi(a^{11},\,a^{12},\ldots,\,a^{33})}{\partial\,a^{i\,j}}\,\mathbf{i}_i\,\mathbf{i}_j(\cdot\,\cdot)(d\,a^{l\,k})\,\mathbf{i}_l\;\mathbf{i}_k,$$

since the expression of the conventional mathematical analysis for the differential of the function of several variables (4.5) follows directly from these relations. As a result, taking into account that $(da^{kl})\,\mathbf{i}_l\mathbf{i}_k=d\mathbf{A}^T$ and $(da^{lk})\,\mathbf{i}_l\mathbf{i}_k=d\mathbf{A}$, we get

---

*If $A^{ij}=C^{ijkl}B_{lk}$ and $A^{ij}$, $B_{lk}$ are the components of the second-rank tensors $\mathbf{A}$ and $\mathbf{B}$, then $C^{ijkl}$ are the components of the fourth-rank tensor $\mathbf{C}$. Special cases: (1) if $a^i=C^{ij}b_j$, where $a^i$ and $b_j$ are the components of the vectors $\mathbf{a}$ and $\mathbf{b}$, then $C^{ij}$ are the components of the second-rank tensor; (2) if $a=c^ib_i$, where $a$ is a scalar and $b_i$ are the components of the vector, then $c^i$ are the components of the vector; (3) if $a=C^{ij}B_{ji}$, where $a$ is a scalar, $B_{ji}$ are the components of the second-rank tensor, then $C^{ij}$ are also the components of the second-rank tensor.

$$d\Phi(\mathbf{A}) = \frac{\partial\Phi(\mathbf{A})}{\partial\mathbf{A}}(\cdot\cdot)d\mathbf{A} = \frac{\partial\Phi(\mathbf{A})}{\partial\mathbf{A}}[\cdot\cdot]d\mathbf{A} = \frac{\partial\Phi(\mathbf{A})}{\partial\mathbf{A}}\cdot\cdot\, d\mathbf{A}^T, \tag{4.6}$$

which allows for the fact that when using the double scalar product of the second rank tensors, the operations $(\cdot\cdot)$ and $[\cdot\cdot]$ are equivalent (see remark after expression (2.12)), $(\cdot\cdot) = \cdot\cdot\,\mathbf{C}_{II}^{IV}\cdot\cdot$ (see remark before expression (2.2)$_1$) and $\mathbf{C}_{II}^{IV}\cdot\cdot\mathbf{A} = \mathbf{A}^T$, where $\mathbf{A}$ is the second-rank tensor (see Table 1). With account of the properties (2.3), expressions (4.6) can be written as

$$d\Phi(\mathbf{A}) = \left(\frac{\partial\Phi(\mathbf{A})}{\partial\mathbf{A}}\right)^T\cdot\cdot\, d\mathbf{A} = \left(\frac{\partial\Phi(\mathbf{A})}{\partial\mathbf{A}}\right)^T(\cdot\cdot)\, d\mathbf{A}^T, \tag{4.7}$$

where

$$\frac{\partial\Phi(\mathbf{A})}{\partial\mathbf{A}} = \frac{\partial\Phi(a^{11}, a^{12},\ldots, a^{33})}{\partial a^{ij}}\,\mathbf{i}_i\,\mathbf{i}_j. \tag{4.8}$$

Knowing tensor (4.8) in the Cartesian basis, we can readily recalculate it in any other basis sets.

Five equivalent expressions (4.6) and (4.7) determine one and the same quantity $d\Phi(\mathbf{A})$, which is the increment of the scalar function of the tensor argument. Moreover, they introduce the other quantity appearing in the right-hand part of these relations in front of the increment of the argument $d\mathbf{A}$. This quantity is a second-rank tensor and is called a derivative of the scalar function of the tensor argument with respect to an argument. There are two such tensor objects, which are coupled by the transpose operation and specified by the scheme of multiplication of the basis vectors of the second-rank tensors, which governs the transition from expression (4.5) to relations (4.6) and (4.7). In mechanics, it is generally agreed that the derivative of the scalar function of the tensor argument with respect to an argument should be defined by the last expression of relation (4.6). However, due to the fact that in this chain of equations the derivative itself remains unchanged, any of the expressions of this relation can be used to determine the desired derivative.

*Remark.* In the literature, it is general practice to introduce the derivative of the scalar function of the tensor argument (as well as for the tensor function of the tensor argument) with respect to an argument with the aid of the Gato derivative (directional derivative)

$$\Delta\Phi(\mathbf{A}) = \frac{d}{ds}\Phi(\mathbf{A} + s\,\Delta\mathbf{A})\,|_{s=0} = \frac{d}{ds}\Phi(a^{ij} + s\,\Delta a^{ij})\,|_{s=0} = \frac{\partial\Phi(a^{ij})}{\partial a^{kp}}\Delta a^{kp} = \Phi_{,\mathbf{A}}(\mathbf{A}) * \Delta\mathbf{A}, \tag{4.9}$$

where $\Delta$ means a small but finite increment of the corresponding quantity, and $*$ is, in compliance with (2.1), either the operation $\cdot\cdot$ or the operation $(\cdot\cdot)$. In the foregoing discussion when deriving the expression for the derivative we immediately proceeded from the next-to-last expression in (4.9).

Hence, expressions (4.3) and (4.6) introduce variations (increments) of the scalar function of the vector and tensor arguments, respectively, whereas relations (4.4) and (4.8) introduce the derivatives of the scalar function of the vector and tensor arguments with respect to this vector or tensor. Note that in many cases in order to find derivatives we need to use only relations (4.3) and (4.6).

Using the last relation in (4.6), we will find the derivatives of the principal invariants of the second rank tensor $\mathbf{A}$ with respect to an argument.

By definition, the first principal invariant of the tensor $\mathbf{A}$ is $I_1(\mathbf{A}) = \mathbf{I}\cdot\cdot\mathbf{A}$. Then, in view of the fact that the unit tensor $\mathbf{I}$ as an object does not change neither in the magnitude nor in the direction due to the invariance of its dyad vectors in the Cartesian basis (see (3.5)), i.e. $d\mathbf{I} = 0$, we obtain

$$dI_1(\mathbf{A}) = d(\mathbf{I}\cdot\cdot\,\mathbf{A}) = d\mathbf{I}\cdot\cdot\,\mathbf{A} + \mathbf{I}\cdot\cdot\, d\mathbf{A} = \mathbf{I}\cdot\cdot\, d\mathbf{A} = \mathbf{I}\cdot\cdot\, d\mathbf{A}^T. \tag{4.10}$$

Here the symmetry of the tensor $\mathbf{I}\,(\mathbf{I} = \mathbf{I}^T)$ and the properties of the tensors described by relation (2.3) are taken into account. Comparing (4.10) and (4.6) we conclude that

$$\frac{\partial I_1(\mathbf{A})}{\partial\mathbf{A}} = \mathbf{I}, \tag{4.11}$$

i.e. the derivative of the first invariant with respect to an argument is equal to a unit tensor.

In the following, we will be in need of derivatives of the first invariant of tensors $\mathbf{A}^2$ and $\mathbf{A}^3$.

$$I_1(\mathbf{A}^2) = \mathbf{I}\cdot\cdot\,\mathbf{A}^2 = \mathbf{I}\cdot\cdot\,(\mathbf{A}\cdot\mathbf{A}) = (\mathbf{I}\cdot\mathbf{A})\cdot\cdot\mathbf{A} = \mathbf{A}\cdot\cdot\mathbf{A}$$

(here the properties (2.4) are used). Then

$$d\,I_1(\mathbf{A}^2) = d\,\mathbf{A}\cdot\cdot\mathbf{A} + \mathbf{A}\cdot\cdot\; d\,\mathbf{A} = 2\,\mathbf{A}\cdot\cdot d\,\mathbf{A} = 2\,\mathbf{A}^T\cdot\cdot d\,\mathbf{A}^T.$$

Comparing this expression with (4.6) we draw a conclusion that

$$\frac{\partial I_1(\mathbf{A}^2)}{\partial \mathbf{A}} = 2\,\mathbf{A}^T. \tag{4.12}$$

Similarly, $I_1(\mathbf{A}^3) = \mathbf{I}\cdot\cdot\mathbf{A}^3 = \mathbf{A}\cdot\cdot\mathbf{A}^2,$

$$d\,I_1(\mathbf{A}^3) = d\,\mathbf{A}\cdot\cdot\mathbf{A}^2 + \mathbf{A}\cdot\cdot\; d\,\mathbf{A}^2 = (\mathbf{A}^2)^T\cdot\cdot\; d\,\mathbf{A}^T + \mathbf{A}\cdot\cdot(d\,\mathbf{A}\cdot\mathbf{A} + \mathbf{A}\cdot d\,\mathbf{A}) = 3\,(\mathbf{A}^2)^T\cdot\cdot d\,\mathbf{A}^T$$

from which we conclude that

$$\frac{\partial I_1(\mathbf{A}^3)}{\partial \mathbf{A}} = 3\,(\mathbf{A}^2)^T. \tag{4.13}$$

Relations (4.11) – (4.13) allow us to write the following expression for an integer positive $n$:

$$\frac{\partial I_1(\mathbf{A}^n)}{\partial \mathbf{A}} = n\;(\mathbf{A}^{n-1})^T.$$

By definition, the second principal invariant of the tensor $\mathbf{A}$ is $I_2(\mathbf{A}) = (I_1^2(\mathbf{A}) - \mathrm{I}_1(\mathbf{A}^2))/2.$ Using relations (4.11) and (4.12), we obtain

$$\frac{\partial I_2(\mathbf{A})}{\partial \mathbf{A}} = I_1(\mathbf{A})\;\mathbf{I} - \mathbf{A}^T. \tag{4.14}$$

Finally, based on the Hamilton-Cayley equality

$$\mathbf{A}^3 - I_1(\mathbf{A})\mathbf{A}^2 + I_2(\mathbf{A})\mathbf{A} - I_3(\mathbf{A})\mathbf{I} = \mathbf{0}, \tag{4.15}$$

and multiplying it twice by the scalar $\mathbf{I}$ (no matter from which side) and taking into consideration that $\mathbf{I}\cdot\cdot\mathbf{I} = I_1(\mathbf{I}) = 3,$ we find the third principal invariant

$$I_3(\mathbf{A}) = \frac{1}{3}\,(I_1(\mathbf{A}^3) - I_1(\mathbf{A})\;I_1(\mathbf{A}^2) + I_2(\mathbf{A})\;I_1(\mathbf{A})).$$

Using expressions (4.11) – (4.13), we have

$$\frac{\partial I_3(\mathbf{A})}{\partial \mathbf{A}} = (\mathbf{A}^2)^T - I_1(\mathbf{A})\;\mathbf{A}^T + I_2(\mathbf{A})\;\mathbf{I}. \tag{4.16}$$

Transposing the Hamilton-Cayley relation (4.15) and taking into account that transposition of the sum of tensors is equal to the sum of the transposed tensors, we obtain

$$(\mathbf{A}^3)^T - I_1(\mathbf{A})(\mathbf{A}^2)^T + I_2(\mathbf{A})\mathbf{A}^T - I_3(\mathbf{A})\mathbf{I} = \mathbf{0}.$$

By performing the scalar multiplication of this equality by $\mathbf{A}^{-T}$, we obtain

$$(\mathbf{A}^2)^T - I_1(\mathbf{A})\mathbf{A}^T + I_2(\mathbf{A})\mathbf{I} = I_3(\mathbf{A})\mathbf{A}^{-T}$$

and then expression (4.16) can be written in a compact form

$$\frac{\partial I_3(\mathbf{A})}{\partial \mathbf{A}} = I_3(\mathbf{A})\;\mathbf{A}^{-T}. \tag{4.17}$$

Thus, expressions (4.11), (4.14) and (4.17) define the derivatives of all three invariants of the second-rank tensor with respect to this tensor. Note that similar expressions are given in [13-18].

Let $\Phi = \Phi(\mathbf{A})$ and $\mathbf{A},$ in its turn be a function of $\mathbf{S}, \mathbf{A} = \mathbf{A}(\mathbf{S}).$ We need to find the derivative $\partial\Phi(\mathbf{A})/\partial\mathbf{S}.$ In accordance with (4.6)

$$d\,\Phi(\mathbf{A}) = \frac{\partial\Phi(\mathbf{A})}{\partial\mathbf{S}}\cdot\cdot\partial\mathbf{S}^T. \tag{4.18}$$

On the other hand,

$$d\,\Phi(\mathbf{A}) = \frac{\partial\Phi(\mathbf{A})}{\partial\mathbf{A}}\cdot\cdot\partial\mathbf{A}^T = \frac{\partial\Phi(\mathbf{A})}{\partial\mathbf{A}}\cdot\cdot(\mathbf{C}_{II}^{IV}\cdot\cdot d\,\mathbf{A}) = (\frac{\partial\Phi(\mathbf{A})}{\partial\mathbf{A}}\cdot\cdot\mathbf{C}_{II}^{IV}\cdot\cdot\frac{\partial\mathbf{A}}{\partial\mathbf{S}})\cdot\cdot\partial\mathbf{S}^T.$$

Here the properties of the tensor $\mathbf{C}_{II}^{IV}$ in the double scalar product (see the third row of Table 1) and relation (5.10) of the next section are taken into account. Comparing the last expression with (4.18) we conclude that

$$\frac{\partial\Phi(\mathbf{A})}{\partial\mathbf{S}} = \frac{\partial\Phi(\mathbf{A})}{\partial\mathbf{A}} \cdot\cdot\, \mathbf{C}_{II}^{IV} \cdot\cdot \frac{\partial\mathbf{A}}{\partial\mathbf{S}}. \tag{4.19}$$

Relation (4.19) can be interpreted as a generalization of the rule of complex function differentiation generally accepted in the mathematical analysis.

## V. Differentiation of Tensor Functions of a Tensor Argument

Now let $\mathbf{\Phi}$ be a tensor function (of second rank) of the tensor argument $\mathbf{A}$ (also of second rank), $\mathbf{\Phi} = \mathbf{\Phi}(\mathbf{A}),$ i.e. each coordinate component of the tensor $\mathbf{\Phi}$ is a function of nine (in general) coordinate components of the tensor $\mathbf{A}.$ Then the expression

$$d\mathbf{\Phi}(\mathbf{A}) = \frac{\partial\Phi^{ij}(a^{11}, a^{12}, \ldots, a^{33})\mathbf{i}_i\mathbf{i}_j}{\partial a^{kp}}\, da^{kp}$$

is a generalization of expression (4.5). Using a special case of the inverse tensor attribute theorem [20] (see footnote on the page with expression (4.1)), the above relation can be represented in the three forms according to schemes (2.5) – (2.7) or (2.9) – (2.11):

$$d\mathbf{\Phi}(\mathbf{A}) = \frac{\partial\Phi^{ij}(a^{11}, a^{12}, \ldots, a^{33})}{\partial a^{kp}}\mathbf{i}_i\mathbf{i}_j\mathbf{i}_k\mathbf{i}_p \cdot\cdot (da^{mn})\mathbf{i}_n\mathbf{i}_m \;\;\Rightarrow\;\; d\mathbf{\Phi}(\mathbf{A}) = \mathbf{L}_I^{IV} \cdot\cdot\, d\mathbf{A}^T, \tag{5.1$_1$}$$

$$d\mathbf{\Phi}(\mathbf{A}) = \frac{\partial\Phi^{ij}(a^{11}, a^{12}, \ldots, a^{33})}{\partial a^{kp}}\mathbf{i}_i\mathbf{i}_j\mathbf{i}_k\mathbf{i}_p (\cdot\cdot)\, (da^{mn})\mathbf{i}_m\mathbf{i}_n \;\;\Rightarrow\;\; d\mathbf{\Phi}(\mathbf{A}) = \mathbf{L}_I^{IV} (\cdot\cdot)\, d\mathbf{A}, \tag{5.1$_2$}$$

$$d\mathbf{\Phi}(\mathbf{A}) = \frac{\partial\Phi^{ij}(a^{11}, a^{12}, \ldots, a^{33})}{\partial a^{kp}}\mathbf{i}_i\mathbf{i}_k\mathbf{i}_p\mathbf{i}_j [\cdot\cdot](da^{mn})\mathbf{i}_m\mathbf{i}_n \;\;\Rightarrow\;\; d\mathbf{\Phi}(\mathbf{A}) = \mathbf{L}_{II}^{IV} [\cdot\cdot]\, d\mathbf{A}. \tag{5.1$_3$}$$

Here $\mathbf{L}_I^{IV}$ is the fourth-rank tensor, representing a derivative $\partial\mathbf{\Phi}(\mathbf{A})/\partial\mathbf{A},$ which can be written in any basis set (in the Cartesian basis it is represented in relations (5.1)$_1$ and (5.1)$_2$). In particular, in the main basis of the initial state it is expressed as $\mathbf{L}_I^{IV} = L^{ijkl}\mathbf{r}_i\mathbf{r}_j\mathbf{r}_k\mathbf{r}_l.$ The tensor $\mathbf{L}_{II}^{IV}$ represents another derivative, which, in contrast to $\mathbf{L}_I^{IV}$, we shall be denote as $(\partial\mathbf{\Phi}(\mathbf{A})/\partial\mathbf{A})^*.$ This tensor is written in the main basis of the initial state as $\mathbf{L}_{II}^{IV} = L^{ijkl}\mathbf{r}_i\mathbf{r}_k\mathbf{r}_l\mathbf{r}_j$ and is the result of a sequential action of the two of transpose operations introduced in [14] on the tensor $\mathbf{L}_I^{IV} : \mathbf{L}_{II}^{IV} = ((\mathbf{L}_I^{IV})^{ti})^{dr},$ where for any tensor of the fourth rank $\mathbf{M}^{IV}$

$$\underset{1234}{(\mathbf{M}^{IV})^{ti}} = \underset{1324}{\mathbf{M}^{IV}}, \qquad \underset{1234}{(\mathbf{M}^{IV})^{dr}} = \underset{1243}{\mathbf{M}^{IV}}, \tag{5.2}$$

In view of the properties of the tensor $\mathbf{C}_{II}^{IV}$ (see Table 1, $\mathbf{A} \cdot\cdot\, \mathbf{C}_{II}^{IV} = \mathbf{C}_{II}^{IV} \cdot\cdot\, \mathbf{A} = \mathbf{A}^T$ for any second-rank tensor $\mathbf{A}$ ), the operation of transposition "$ti$" can be represented as

$$(\mathbf{M}^{IV})^{ti} = \mathbf{C}_{II}^{IV} \overset{2}{\odot} \mathbf{M}^{IV} \qquad or \qquad (\mathbf{M}^{IV})^{ti} = \mathbf{M}^{IV} \overset{3}{\circledast} \mathbf{C}_{II}^{IV}. \tag{5.3}$$

Here, $\overset{n}{\odot}$ is the operation of the positional double scalar multiplication $\cdot\cdot$ of the tensor, which is to the left in this operation, by the $n$-th and $n+1$-th basis vectors of the tensor, which is to the right in this operation:

$$\mathbf{C}_{II}^{IV} \overset{2}{\odot} \mathbf{M}^{IV} = M^{ijkl}\mathbf{r}_i(\mathbf{C}_{II}^{IV} \cdot\cdot\, \mathbf{r}_j\mathbf{r}_k)\mathbf{r}_l = M^{ijkl}\mathbf{r}_i\mathbf{r}_k\mathbf{r}_j\mathbf{r}_l.$$

$\overset{n}{\circledast}$ is the operation of the positional double product $\cdot\cdot$ of the tensor, which is to the right in this operation, by the $n$-th and $n-1$-th basis vectors of the tensor, which is to the left in this operation:

$$\mathbf{M}^{IV} \overset{3}{\circledast} \mathbf{C}_{II}^{IV} = M^{ijkl}\mathbf{r}_i(\mathbf{r}_j\mathbf{r}_k \cdot\cdot\, \mathbf{C}_{II}^{IV})\mathbf{r}_l = M^{ijkl}\mathbf{r}_i\mathbf{r}_k\mathbf{r}_j\mathbf{r}_l.$$

Using the introduced operators the transposition operation "$dr$" is represented as

$$(\mathbf{M}^{IV})^{dr} = \mathbf{C}_{II}^{IV} \overset{3}{\odot} \mathbf{M}^{IV} \qquad or \qquad (\mathbf{M}^{IV})^{dr} = \mathbf{M}^{IV} \cdot\cdot\, \mathbf{C}_{II}^{IV} \tag{5.4}$$

and with account of representations (5.3) and (5.4), we finally arrive at the following four equivalent expressions for the tensor $\mathbf{L}_{II}^{IV}$

$$\mathbf{L}_{II}^{IV} = ((\mathbf{L}_{I}^{IV})^{ti})^{dr} = \mathbf{C}_{II}^{IV} \overset{3}{\odot} (\mathbf{C}_{II}^{IV} \overset{2}{\odot} \mathbf{L}_{I}^{IV}) = \mathbf{C}_{II}^{IV} \overset{3}{\odot} (\mathbf{L}_{I}^{IV} \overset{3}{\circledast} \mathbf{C}_{II}^{IV}) =$$

$$= (\mathbf{L}_{I}^{IV} \overset{3}{\circledast} \mathbf{C}_{II}^{IV}) \cdot\cdot \mathbf{C}_{II}^{IV} = (\mathbf{C}_{II}^{IV} \overset{2}{\odot} \mathbf{L}_{I}^{IV}) \cdot\cdot \mathbf{C}_{II}^{IV}. \tag{5.5}$$

Note that the equality of derivatives in relations $(5.1)_1$ and $(5.1)_2$, follows from the previously determined relationship between the operators $(\cdot\cdot) = \cdot\cdot\, \mathbf{C}_{II}^{IV} \,\cdot\cdot$ (see the remark before expression $(2.2)_1$).
Below we will be used expression $(5.1)_2$ in the form

$$d\mathbf{\Phi}(\mathbf{A}) = \frac{\partial \mathbf{\Phi}(\mathbf{A})}{\partial \mathbf{A}} \cdot\cdot\, d\mathbf{A}^T \tag{5.6}$$

for the construction of derivatives and development of the rules for their derivation. However, before proceeding to their consideration it will be necessary to obtain a number of relationships required for our further discussion. From the chain of equalities

$$\mathbf{A}[\cdot\cdot]\mathbf{L}_{II}^{IV} = A^{ij}\mathbf{r}_i \cdot \mathbf{L}_{II}^{IV} \cdot \mathbf{r}_j = A^{ij}\mathbf{r}_i \cdot L_{mnkp}\, \mathbf{r}^m\mathbf{r}^k\mathbf{r}^p\mathbf{r}^n \cdot \mathbf{r}_j = A^{ij}\, L_{ijkp}\, \mathbf{r}^k\mathbf{r}^p = \mathbf{A}(\cdot\cdot)\mathbf{L}_{I}^{IV} = \mathbf{A} \cdot\cdot\, \mathbf{C}_{II}^{IV} \cdot\cdot\, \mathbf{L}_{I}^{IV},$$

where $\mathbf{A}$ is any second-rank tensor, we obtain the following relation

$$\mathbf{A}[\cdot\cdot]\mathbf{L}_{II}^{IV} = \mathbf{A}(\cdot\cdot)\mathbf{L}_{I}^{IV} = \mathbf{A} \cdot\cdot\, \mathbf{C}_{II}^{IV} \cdot\cdot\, \mathbf{L}_{I}^{IV}, \tag{5.7}$$

and the chain of equalities

$${}^{(1)}\mathbf{L}_{II}^{IV}[\cdot\cdot]\,{}^{(2)}\mathbf{L}_{II}^{IV} = {}^{(1)}L^{ijkl}\mathbf{r}_i(\mathbf{r}_k \cdot {}^{(2)}\mathbf{L}_{II}^{IV} \cdot \mathbf{r}_l)\mathbf{r}_j = {}^{(1)}L^{ijkl}\mathbf{r}_i(\mathbf{r}_k \cdot {}^{(2)}L_{mnsq}\mathbf{r}^m\mathbf{r}^s\mathbf{r}^q\mathbf{r}^n \cdot \mathbf{r}_l)\mathbf{r}_j = {}^{(1)}L^{ijkl}\,{}^{(2)}L_{klsq}\,\mathbf{r}_i\mathbf{r}^s\mathbf{r}^q\mathbf{r}_j,$$

where ${}^{(1)}\mathbf{L}_{II}^{IV}$ and ${}^{(2)}\mathbf{L}_{II}^{IV}$ are the two tensors $\mathbf{L}_{II}^{IV}$ (different in the general case), leads to the relation

$${}^{(1)}\mathbf{L}_{II}^{IV}[\cdot\cdot]\,{}^{(2)}\mathbf{L}_{II}^{IV} = [({}^{(1)}\mathbf{L}_{I}^{IV}(\cdot\cdot)\,{}^{(2)}\mathbf{L}_{I}^{IV})^{ti}]^{dr}. \tag{5.8}$$

To perform transpose operations in the last equality, schemes (5.2) or relations (5.3) – (5.5) can be used.

Based on relation (5.6) we will obtain some useful expressions. Let $\mathbf{\Phi} = \mathbf{\Phi}(\mathbf{A})$ and $\mathbf{A}$, in turn, be a function of $\mathbf{S}$, $\mathbf{A} = \mathbf{A}(\mathbf{S})$. So, it is necessary to find the derivative $\partial \mathbf{\Phi}(\mathbf{A}) / \partial \mathbf{S}$. According to (5.6) we have

$$d\mathbf{\Phi}(\mathbf{A}) = \frac{\partial \mathbf{\Phi}(\mathbf{A})}{\partial \mathbf{S}} \cdot\cdot\, d\mathbf{S}^T. \tag{5.9}$$

On the other hand,

$$d\mathbf{\Phi}(\mathbf{A}) = \frac{\partial \mathbf{\Phi}(\mathbf{A})}{\partial \mathbf{A}} \cdot\cdot\, d\mathbf{A}^T = \frac{\partial \mathbf{\Phi}(\mathbf{A})}{\partial \mathbf{A}} \cdot\cdot\, d(\mathbf{C}_{II}^{IV} \cdot\cdot\, \mathbf{A}) =$$

$$= \frac{\partial \mathbf{\Phi}(\mathbf{A})}{\partial \mathbf{A}} \cdot\cdot\, \mathbf{C}_{II}^{IV} \cdot\cdot\, d\mathbf{A} = \left(\frac{\partial \mathbf{\Phi}(\mathbf{A})}{\partial \mathbf{A}} \cdot\cdot\, \mathbf{C}_{II}^{IV} \cdot\cdot\, \frac{\partial \mathbf{A}}{\partial \mathbf{S}}\right) \cdot\cdot\, d\mathbf{S}^T.$$

Here we take into account the relation (5.6) $d\mathbf{A}(\mathbf{S}) = \partial \mathbf{A}(\mathbf{S}) / \partial \mathbf{S} \cdot\cdot\, d\mathbf{S}^T$ and the properties of the tensor $\mathbf{C}_{II}^{IV}$ given in Table 1, its invariability either in the magnitude or direction due to invariability of its dyad vectors in the Cartesian basis (see (3.6)), i.e. $d\mathbf{C}_{II}^{IV} = 0$ from which it follows that $d(\mathbf{C}_{II}^{IV} \cdot\cdot\, \mathbf{A}) = \mathbf{C}_{II}^{IV} \cdot\cdot\, d\mathbf{A}$. Comparing the last expression in the above chain of equalities with (5.9), we conclude that

$$\frac{\partial \mathbf{\Phi}(\mathbf{A})}{\partial \mathbf{S}} = \frac{\partial \mathbf{\Phi}(\mathbf{A})}{\partial \mathbf{A}} \cdot\cdot\, \mathbf{C}_{II}^{IV} \cdot\cdot\, \frac{\partial \mathbf{A}}{\partial \mathbf{S}}. \tag{5.10}$$

Expression (5.10) can be interpreted as a generalization of the rule of complex function differentiation widely used in the mathematical analysis.

Now let $\mathbf{\Phi}(\mathbf{S}) = \mathbf{A}(\mathbf{S}) \cdot \mathbf{B}(\mathbf{S})$. Then, on the one hand, according to (5.6)

$$d\mathbf{\Phi}(\mathbf{S}) = \frac{\partial \mathbf{\Phi}(\mathbf{S})}{\partial \mathbf{S}} \cdot\cdot\, d\mathbf{S}^T, \tag{5.11}$$

and on the other hand

$$d\mathbf{\Phi}(\mathbf{S}) = d\mathbf{A}\cdot\mathbf{B} + \mathbf{A}\cdot d\mathbf{B} = (\frac{\partial\mathbf{A}}{\partial\mathbf{S}}\cdot\cdot\, d\mathbf{S}^T)\cdot\mathbf{B} + \mathbf{A}\cdot(\frac{\partial\mathbf{B}}{\partial\mathbf{S}}\cdot\cdot\, d\mathbf{S}^T), \tag{5.12}$$

where the expressions of the type (5.6) are used for $d\mathbf{A}(\mathbf{S})$ and $d\mathbf{B}(\mathbf{S})$. Here, we need to introduce a scheme of scalar products of the basis vectors of the last two summands, taking into account that derivatives entering into these summands are the fourth-rank tensors (with four basis vectors). The left-hand part of relation (5.13) given below contains a scheme of relationship (scalar product) between the basis vectors of tensors entering into the last summand in (5.12). The horizontal brackets indicate for which basis vectors and of which tensors the scalar-wise multiplication is performs. For example, the second basis vector of the second-rank tensor $\mathbf{A}$ is multiplied by the first basis vector of the fourth-rank tensor $\partial\mathbf{B}/\partial\mathbf{S}$, whereas the third and the fourth basis vectors of this tensor are scalar-wise related to the second and first basis vectors of the tensor $d\mathbf{S}^T$, respectively. The structure of relation between the basis vectors of all tensors entering into the right part of relation (5.13) is similar to the structure of the left-hand part. Therefore, the sign of equality is entered between them. As a result, we obtain a tensor of the second rank, whose first basis vector is the first basis vector of the tensor $\mathbf{A}$, and the second basis vector is the second basis vector of the fourth-rank tensor $\partial\mathbf{B}/\partial\mathbf{S}$.

$$\mathbf{A}\cdot(\frac{\partial\mathbf{B}}{\partial\mathbf{S}}\cdot\cdot\, d\mathbf{S}^T) = (\mathbf{A}\cdot\frac{\partial\mathbf{B}}{\partial\mathbf{S}})\cdot\cdot\, d\mathbf{S}^T. \tag{5.13}$$

The structure of the relationship between the basis vectors of the first term in the right-hand part of relation (5.12) is shown below:

$$(\frac{\partial\mathbf{A}}{\partial\mathbf{S}}\cdot\cdot\, d\mathbf{S}^T)\cdot\mathbf{B}. \tag{5.14}$$

As a result, we obtain the second-rank tensor, the first basis vector of which is the first basis vector of the fourth-rank tensor $\partial\mathbf{A}/\partial\mathbf{S}$ and the second basis vector is the second basis vector of the tensor $\mathbf{B}$. Earlier we introduced the operation of positional double scalar product. Now we use the operation of simple positional scalar multiplication $\overset{n}{*}$ introduced in [19]. Let $\mathbf{H}^{IV}$ be the tensor of the fourth rank and $\mathbf{D}$ is the tensor of the second rank. Then the operation $\mathbf{H}^{IV}\overset{n}{*}\mathbf{D}$ means the scalar multiplication of the $n$-th basis vector of the fourth-rank tensor $\mathbf{H}^{IV}$ on the left by the second-rank tensor $\mathbf{D}$. For example, if the tensor $\mathbf{H}^{IV}$ is written as $\mathbf{H}^{IV} = H^{ijkl}\mathbf{r}_i\mathbf{r}_j\mathbf{r}_k\mathbf{r}_l$ and tensor $\mathbf{D}$ is written as $\mathbf{D} = D^{mn}\mathbf{r}_m\mathbf{r}_n$, then

$$\mathbf{H}^{IV}\overset{2}{*}\mathbf{D} = H^{ijkl}\mathbf{r}_i(\mathbf{r}_j\cdot\mathbf{D})\mathbf{r}_k\mathbf{r}_l = H^{ijkl}\mathbf{r}_i(\mathbf{r}_j\cdot D^{mn}\mathbf{r}_m\mathbf{r}_n)\mathbf{r}_k\mathbf{r}_l = H^{ijkl}D^{mn}g_{jm}\mathbf{r}_i\mathbf{r}_n\mathbf{r}_k\mathbf{r}_l.$$

Using this operation, expression (5.14) is rewritten in the following form:

$$(\frac{\partial\mathbf{A}}{\partial\mathbf{S}}\overset{2}{*}\mathbf{B})\cdot\cdot\, d\mathbf{S}^T. \tag{5.15}$$

The scheme of relationship between the basis vectors in this expression is identical to the scheme in (5.14). Therefore, with account of expression (5.13) and (5.15) relation (5.12) can be represented as

$$d\mathbf{\Phi}(\mathbf{S}) = (\frac{\partial\mathbf{A}}{\partial\mathbf{S}}\overset{2}{*}\mathbf{B} + \mathbf{A}\cdot\frac{\partial\mathbf{B}}{\partial\mathbf{S}})\cdot\cdot\, d\mathbf{S}^T.$$

Comparing this expression with (5.11), we may conclude that

$$\frac{\partial(\mathbf{A}(\mathbf{S})\cdot\mathbf{B}(\mathbf{S}))}{\partial\mathbf{S}} = \frac{\partial\mathbf{A}}{\partial\mathbf{S}}\overset{2}{*}\mathbf{B} + \mathbf{A}\cdot\frac{\partial\mathbf{B}}{\partial\mathbf{S}}. \tag{5.16}$$

By its structure, (5.16) resembles, to a certain extent, a conventional rule of the mathematical analysis for product differentiation:

$$\frac{\partial}{\partial x}(y(x)\, f(x)) = \frac{\partial y(x)}{\partial x}\, f(x) + y(x)\frac{\partial f(x)}{\partial x}.$$

In the following, we will consider some simple examples. Let $\mathbf{\Phi}(\mathbf{A}) = \mathbf{A}.$ Then relation (5.6) yields

$$d\,\mathbf{A} = \frac{\partial\,\mathbf{A}}{\partial\,\mathbf{A}} \cdot\cdot\, d\,\mathbf{A}^T = \mathbf{C}_{II}^{IV} \cdot\cdot\, d\,\mathbf{A}^T$$

(see the properties of the tensor $\mathbf{C}_{II}^{IV}$ in Table 1), which suggests that

$$\frac{\partial\,\mathbf{A}}{\partial\,\mathbf{A}} = \mathbf{C}_{II}^{IV}. \tag{5.17}$$

Now, let $\mathbf{\Phi}(\mathbf{A}) = \mathbf{A}^T.$ Then, according to relation (5.6), we have

$$d\,\mathbf{A}^T = \frac{\partial\,\mathbf{A}^T}{\partial\,\mathbf{A}} \cdot\cdot\, d\,\mathbf{A}^T = \mathbf{C}_{III}^{IV} \cdot\cdot\, d\,\mathbf{A}^T$$

(see the properties of the tensor $\mathbf{C}_{III}^{IV}$ in Table 1), which suggests that

$$\frac{\partial\,\mathbf{A}^T}{\partial\,\mathbf{A}} = \mathbf{C}_{III}^{IV}. \tag{5.18}$$

If $\mathbf{\Phi}(\mathbf{A}) = \mathbf{A}^2 = \mathbf{A}\cdot\mathbf{A},$ then expression (5.16), taking into account (5.17), gives the following derivative:

$$\frac{\partial\,\mathbf{A}^2}{\partial\,\mathbf{A}} = \frac{\partial\,(\mathbf{A}\cdot\mathbf{A})}{\partial\,\mathbf{A}} = \frac{\partial\,\mathbf{A}}{\partial\,\mathbf{A}} \overset{2}{*}\ \mathbf{A} + \mathbf{A}\cdot\frac{\partial\,\mathbf{A}}{\partial\,\mathbf{A}} = \mathbf{C}_{II}^{IV} \overset{2}{*}\ \mathbf{A} + \mathbf{A}\cdot\,\mathbf{C}_{II}^{IV}. \tag{5.19}$$

This derivative, again, can be imagined (after replacing $\mathbf{C}_{II}^{IV}$ in (5.19) by $\mathbf{I}$ and the positional multiplication by a simple scalar product) as being consistent with the expression of conventional mathematical analysis $\partial x^2 / \partial x = 2\,x.$

Let $\mathbf{\Phi}(\mathbf{A}) = \mathbf{A}^{-1}.$ Then, on the one hand,

$$\frac{\partial\,(\mathbf{A}\cdot\mathbf{A}^{-1})}{\partial\,\mathbf{A}} = 0,$$

since $\mathbf{A}\cdot\mathbf{A}^{-1} = \mathbf{I}$ and, on the other hand, by virtue of (5.16) and (5.17),

$$\frac{\partial\,\mathbf{A}\cdot\mathbf{A}^{-1}}{\partial\,\mathbf{A}} = \mathbf{C}_{II}^{IV} \overset{2}{*}\ \mathbf{A}^{-1} + \mathbf{A}\cdot\frac{\partial\ \mathbf{A}^{-1}}{\partial\ \mathbf{A}}.$$

From these two expressions we obtain

$$\frac{\partial\,\mathbf{A}^{-1}}{\partial\,\mathbf{A}} = -\,\mathbf{A}^{-1}\cdot\,\mathbf{C}_{II}^{IV} \overset{2}{*}\ \mathbf{A}^{-1}. \tag{5.20}$$

Again, this can be conceived as being consistent with the well-known expression of conventional mathematical analysis $\partial x^{-1} / \partial x = -\,x^{-2}.$

Finally, let $\mathbf{\Phi}(\mathbf{D}) = \Psi(\mathbf{A})\,\mathbf{\Lambda}(\mathbf{B}),\ \mathbf{A} = \mathbf{A}(\mathbf{S}),\ \mathbf{B} = \mathbf{B}(\mathbf{S}),$ i.e. the tensor function of a tensor argument is the product of a scalar function of the tensor argument $\Psi(\mathbf{A})$ and a tensor (second-rank) function of the tensor argument $\mathbf{\Lambda}(\mathbf{B}).$ The arguments, in turn, depend on the tensor $\mathbf{S}.$ We define the derivative $\partial\mathbf{\Phi}(\mathbf{D})/\partial\mathbf{S}$ as

$$d\,\mathbf{\Phi}(\mathbf{D}) = d\,(\Psi(\mathbf{A})\,\mathbf{\Lambda}(\mathbf{B})) = d\,(\Psi(\mathbf{A}))\,\mathbf{\Lambda}(\mathbf{B}) + \Psi(\mathbf{A})\,d\,(\mathbf{\Lambda}(\mathbf{B})) =$$

$$= (\frac{\partial\Psi(\mathbf{A})}{\partial\mathbf{S}} \cdot\cdot\, d\mathbf{S}^T)\ \mathbf{\Lambda}(\mathbf{B}) + \Psi(\mathbf{A})\,(\frac{\partial\mathbf{\Lambda}(\mathbf{B})}{\partial\mathbf{S}} \cdot\cdot\, d\mathbf{S}^T) =$$

$$=(\mathbf{\Lambda}(\mathbf{B})\frac{\partial\Psi(\mathbf{A})}{\partial\mathbf{S}}+\Psi(\mathbf{A})\frac{\partial\mathbf{\Lambda}(\mathbf{B})}{\partial\mathbf{S}})\cdot\cdot\, d\mathbf{S}^{T}. \qquad (5.21)$$

This expression takes into account the fact that the expression in parentheses in the first term of the second row is a scalar quantity. From a comparison of (5.21) with (5.6), it follows that

$$\frac{\partial}{\partial\mathbf{S}}(\Psi(\mathbf{A})\,\mathbf{\Lambda}(\mathbf{B}))=\mathbf{\Lambda}(\mathbf{B})\frac{\partial\Psi(\mathbf{A})}{\partial\mathbf{S}}+\Psi(\mathbf{A})\frac{\partial\mathbf{\Lambda}(\mathbf{B})}{\partial\mathbf{S}}. \qquad (5.22)$$

Note that the sequence order of the multiplicands in the tensor product in the first summand of the right-hand part of the last relation cannot be changed, since the tensor product is noncommutative. Of course, the scalar factor $\Psi(\mathbf{A})$ in the second summand can appear anywhere in this term. Taking into account equalities (4.19) and (5.10), (5.22) can be represented as

$$\frac{\partial}{\partial\mathbf{S}}(\Psi(\mathbf{A})\,\mathbf{\Lambda}(\mathbf{B}))=\mathbf{\Lambda}(\mathbf{B})(\frac{\partial\Psi(\mathbf{A})}{\partial\mathbf{A}}\cdot\cdot\,\mathbf{C}_{II}^{IV}\cdot\cdot\frac{\partial\mathbf{A}}{\partial\mathbf{S}})+\Psi(\mathbf{A})(\frac{\partial\mathbf{\Lambda}(\mathbf{B})}{\partial\mathbf{B}}\cdot\cdot\,\mathbf{C}_{II}^{IV}\cdot\cdot\frac{\partial\mathbf{B}}{\partial\mathbf{S}}). \qquad (5.23)$$

## VI. Rules for Differentiation of Tensor Functions of a Tensor Argument and Derivatives Obtained in Other Publications

As noted above, in [8-10] (in the Introduction these publications are related to the first group) the operation of the double scalar product of tensors $(\cdot\cdot)$ is used (see (2.9), scheme (2.6) and the second scheme in (2.8)), and the derivative of the tensor function with respect to a tensor argument $\mathbf{L}_{I}^{IV}$ is determined by relation $(5.1)_2$. In [11-14] (in the Introduction they are referred to the second group) the operation of double scalar product of tensors $[\cdot\cdot]$ is applied (see (2.11), scheme (2.5) and the first scheme in (2.8)) and the derivative of the tensor function with respect to a tensor argument $\mathbf{L}_{II}^{IV}$ is specified by relation $(5.1)_3$. In this work, the following relationships between the operators $[\cdot\cdot]$, $(\cdot\cdot)$ and $\cdot\cdot$ are established and presented in Table 2.

**Table 2** *The relationships between the operators* $[\cdot\cdot]$, $(\cdot\cdot)$ *and* $\cdot\cdot$

| Number | Linkage | Comment |
|---|---|---|
| 1. | For a double scalar product of the second rank tensors the operations $(\cdot\cdot)$ and $[\cdot\cdot]$ are equivalent one another | Remark after relation (2.12) |
| 2. | $(\cdot\cdot)=\cdot\cdot\,\mathbf{C}_{II}^{IV}\cdot\cdot$ | Remark before relation $(2.2)_1$ |
| 3. | $\mathbf{C}_{II}^{IV}\cdot\cdot\mathbf{A}=\mathbf{A}\cdot\cdot\mathbf{C}_{II}^{IV}=\mathbf{A}^{T}$, where $\mathbf{A}$ is a second-rank tensor | Table 1 |
| 4. | $\mathbf{A}[\cdot\cdot]\mathbf{L}_{II}^{IV}=\mathbf{A}(\cdot\cdot)\mathbf{L}_{I}^{IV}=\mathbf{A}\cdot\cdot\,\mathbf{C}_{II}^{IV}\cdot\cdot\,\mathbf{L}_{I}^{IV}$ for any second-rank tensor $\mathbf{A}$ | Relation (5.7) |
| 5. | ${}^{(1)}\mathbf{L}_{II}^{IV}[\cdot\cdot]\,{}^{(2)}\mathbf{L}_{II}^{IV}=[({}^{(1)}\mathbf{L}_{I}^{IV}(\cdot\cdot)\,{}^{(2)}\mathbf{L}_{I}^{IV})^{ti}]^{dr}$ for any two tensors $\mathbf{L}_{II}^{IV}$ which are different in the general case | Relation (5.8)). The schemes (5.2) or the relations (5.3) – (5.5) can be used to perform transpose operations in the last equality |

Using these relationships, we will try to establish the correspondence between the rules for differentiation of the scalar and tensor functions of the tensor argument obtained in [15-18], as well as in this article (in the Introduction these publications are referred to the third group), and the rules obtained in the publications of the first two groups, and also between the resulting derivatives. As noted earlier, and follows from the links established in paragraphs 1 – 3 of Table 2, the rules for differentiating the scalar functions of the tensor argument and the resulting derivatives are identical in all publications. The correspondence between the rules for differentiation of the tensor functions of the tensor argument and between the derivatives obtained in the publications of these three groups will be varified in compliance with the order of their appearance in this article.

*6.1. Derivative* $\Phi(\mathbf{A})_{,\mathbf{S}}$. The first in this sequence is relation (4.19), which determines the derivative of the scalar function $\Phi$ depending on the second rank tensor $\mathbf{A}$, which, in turn, is a function of $\mathbf{S}$, $\mathbf{A}=\mathbf{A}(\mathbf{S})$: $\Phi(\mathbf{A})_{,\mathbf{S}}=\Phi(\mathbf{A})_{,\mathbf{A}}\cdot\cdot\,\mathbf{C}_{II}^{IV}\cdot\cdot\mathbf{A}_{,\mathbf{S}}$. In the works of the first group this representation corresponds to the expression $\Phi(\mathbf{A})_{,\mathbf{S}}=\Phi(\mathbf{A})_{,\mathbf{A}}(\cdot\cdot)\,\mathbf{A}_{,\mathbf{S}}$ while in the works of the second group – to the expression $\mathbf{\Phi}(\mathbf{A})_{,\mathbf{S}} = \mathbf{\Phi}(\mathbf{A})_{,\mathbf{A}}[\cdot\cdot]\mathbf{A}^*_{,\mathbf{S}}$ (hereinafter, notation $\mathbf{A}_{,\mathbf{S}} = \partial\mathbf{A}/\partial\mathbf{S}$ is used and in this case $\mathbf{A}_{,\mathbf{S}} = \mathbf{L}^{IV}_{I}$, $\mathbf{A}^*_{,\mathbf{S}} = \mathbf{L}^{IV}_{II}$, see definitions of derivatives $(5.1)_1$, $(5.1)_2$ and $(5.1)_3$). The equivalence of all these expressions follows from the row 4 of Table 2.

*6.2. Derivative* $\mathbf{\Phi}(\mathbf{A})_{,\mathbf{S}}$. Relation (5.10) defines the derivative of the second-rank tensor function $\mathbf{\Phi}$ of the tensor argument $\mathbf{A}$ (also of the second rank) which, in turn, is a function of $\mathbf{S}, \mathbf{A} = \mathbf{A}(\mathbf{S})$: $\mathbf{\Phi}(\mathbf{A})_{,\mathbf{S}} = \mathbf{\Phi}(\mathbf{A})_{,\mathbf{A}} \cdot\cdot \mathbf{C}^{IV}_{II} \cdot\cdot \mathbf{A}_{,\mathbf{S}}$. In the first group of publications, this form of representation is consistent with the expression $\mathbf{\Phi}(\mathbf{A})_{,\mathbf{S}} = \mathbf{\Phi}(\mathbf{A})_{,\mathbf{A}}(\cdot\cdot)\mathbf{A}_{,\mathbf{S}}$ which, according to item 2 of Table 2, is fully consistent with the previous one. In the works of the second group relation (5.10) corresponds to the expression $\mathbf{\Phi}(\mathbf{A})^*_{,\mathbf{S}} = \mathbf{\Phi}(\mathbf{A})^*_{,\mathbf{A}}[\cdot\cdot]\mathbf{A}^*_{,\mathbf{S}}$ which with account of item 5 of Table 2 can be written as $\mathbf{\Phi}(\mathbf{A})^*_{,\mathbf{S}} = [(\mathbf{\Phi}(\mathbf{A})_{,\mathbf{A}}(\cdot\cdot)\mathbf{A}_{,\mathbf{S}})^{ti}]^{dr}$ which, in turn, according to row 2 of Table 2, can be written as $\mathbf{\Phi}(\mathbf{A})^*_{,\mathbf{S}} = [(\mathbf{\Phi}(\mathbf{A})_{,\mathbf{A}} \cdot\cdot \mathbf{C}^{IV}_{II} \cdot\cdot \mathbf{A}_{,\mathbf{S}})^{ti}]^{dr}$. These expressions coincide with the relation of the first group and relation (5.10) to within the operations of transposition "$ti$" and "$dr$". This is just the relationship that is established between the derivatives $\mathbf{L}^{IV}_{II}$ and $\mathbf{L}^{IV}_{I}$ and defined by the relations $(5.1)_1$, $(5.1)_2$ and $(5.1)_3$.

*6.3. Derivative* $(\mathbf{A}(\mathbf{S}) \cdot \mathbf{B}(\mathbf{S}))_{,\mathbf{S}}$. The derivative of the scalar product of two tensor (second-rank) functions of the same tensor argument (of the second rank, too) with respect to this argument is determined by relation (5.16): $(\mathbf{A}(\mathbf{S}) \cdot \mathbf{B}(\mathbf{S}))_{,\mathbf{S}} = \mathbf{A}_{,\mathbf{S}} \overset{2}{*} \mathbf{B} + \mathbf{A} \cdot \mathbf{B}_{,\mathbf{S}}$. In the works of the first group this representation corresponds to the expression

$$(\mathbf{A}(\mathbf{S}) \cdot \mathbf{B}(\mathbf{S}))_{,\mathbf{S}} = (\mathbf{A} \boxtimes \mathbf{I})(\cdot\cdot)\mathbf{B}_{,\mathbf{S}} + (\mathbf{I} \boxtimes \mathbf{B}^T)(\cdot\cdot)\mathbf{A}_{,\mathbf{S}}, \tag{6.1}$$

which, taking into account item 2 of Table 2, can be rewritten as

$$(\mathbf{A}(\mathbf{S}) \cdot \mathbf{B}(\mathbf{S}))_{,\mathbf{S}} = (\mathbf{A} \boxtimes \mathbf{I}) \cdot\cdot \mathbf{C}^{IV}_{II} \cdot\cdot \mathbf{B}_{,\mathbf{S}} + (\mathbf{I} \boxtimes \mathbf{B}^T) \cdot\cdot \mathbf{C}^{IV}_{II} \cdot\cdot \mathbf{A}_{,\mathbf{S}}. \tag{6.1$_1$}$$

Let us expand the first term of the last expression. Taking into account that $\mathbf{A} \boxtimes \mathbf{I} = A^{ij}\mathbf{r}_i\mathbf{r}^k\mathbf{r}_j\mathbf{r}_k$, $\mathbf{C}^{IV}_{II} = \mathbf{r}^i\mathbf{r}^j\mathbf{r}_i\mathbf{r}_j$, $\mathbf{B}_{,\mathbf{S}} = B_{ijkl}\mathbf{r}^i\mathbf{r}^j\mathbf{r}^k\mathbf{r}^l$, we obtain $(\mathbf{A} \boxtimes \mathbf{I}) \cdot\cdot \mathbf{C}^{IV}_{II} \cdot\cdot \mathbf{B}_{,\mathbf{S}} = A^{ij}B_{jkmn}\mathbf{r}_i\mathbf{r}^k\mathbf{r}^m\mathbf{r}^n = \mathbf{A} \cdot \mathbf{B}_{,\mathbf{S}}$, which corresponds to the second term of expression (5.16). Given that $\mathbf{I} \boxtimes \mathbf{B}^T = B^{ij}\mathbf{r}^k\mathbf{r}_j\mathbf{r}_k\mathbf{r}_i$ and representing $\mathbf{A}_{,\mathbf{S}}$ similarly to $\mathbf{B}_{,\mathbf{S}}$ in the previous case, we have $(\mathbf{I} \boxtimes \mathbf{B}^T) \cdot\cdot \mathbf{C}^{IV}_{II} \cdot\cdot \mathbf{A}_{,\mathbf{S}} = A_{kimn}B^{ij}\mathbf{r}^k\mathbf{r}_j\mathbf{r}^m\mathbf{r}^n$. Using the operation of positional scalar multiplication introduced after relation (5.14), we can represent the last expression as $\mathbf{A}_{,\mathbf{S}} \overset{2}{*} \mathbf{B}$, which corresponds to the first term in (5.16). This confirms a complete correspondence between relation (5.16) and the relevant expression (6.1) from the works of the first group.

In the works of the second group the derivative of the scalar product of two second-rank tensor functions of the same tensor argument of the second rank, too, with respect to this argument is written as

$$(\mathbf{A}(\mathbf{S}) \cdot \mathbf{B}(\mathbf{S}))^*_{,\mathbf{S}} = \mathbf{A}^*_{,\mathbf{S}} \cdot \mathbf{B} + \mathbf{A} \cdot \mathbf{B}^*_{,\mathbf{S}}. \tag{6.2}$$

In the coordinate representation $\mathbf{A}^*_{,\mathbf{S}} = A^{ijkp}\mathbf{r}_i\mathbf{r}_k\mathbf{r}_p\mathbf{r}_j$ (the tensor $\mathbf{B}^*_{,\mathbf{S}}$ has a similar representation) and $\mathbf{B} = B_{mn}\mathbf{r}^m\mathbf{r}^n$ (the tensor $\mathbf{A}$ has a similar representation). As a result, $\mathbf{A}^*_{,\mathbf{S}} \cdot \mathbf{B} = A^{ijkp}B_{jn}\mathbf{r}_i\mathbf{r}_k\mathbf{r}_p\mathbf{r}^n$, which can be written as $[(\mathbf{A}_{,\mathbf{S}} \overset{2}{*} \mathbf{B})^{ti}]^{dr}$, and $\mathbf{A} \cdot \mathbf{B}^*_{,\mathbf{S}} = A_{mi}B^{ijkp}\mathbf{r}^m\mathbf{r}_k\mathbf{r}_p\mathbf{r}_j$, which can be written as $[(\mathbf{A} \cdot \mathbf{B}_{,\mathbf{S}})^{ti}]^{dr}$. These expressions coincide with the corresponding summands of the first group and with the summands of expression (5.16) to within the transpose operations "$ti$" and "$dr$". With account of the fact that the transpose of a sum is the sum of the transposed variables, we obtain $[(\mathbf{A}_{,\mathbf{S}} \overset{2}{*} \mathbf{B} + \mathbf{A} \cdot \mathbf{B}_{,\mathbf{S}})^{ti}]^{dr}$, which suggests that derivative $\mathbf{L}^{IV}_{II}$ should be defined in terms of derivative $\mathbf{L}^{IV}_{I}$ by relations $(5.1)_1$, $(5.1)_2$ and $(5.1)_3$.

*6.4. Derivatives* $\mathbf{A}_{,\mathbf{A}}$ *and* $\mathbf{A}^T_{,\mathbf{A}}$. The derivative of the second-rank tensor $\mathbf{A}$ with respect to the tensor $\mathbf{A}$ is defined by (5.17): $\mathbf{A}_{,\mathbf{A}} = \mathbf{C}^{IV}_{II}$. In the works of the first group this expression corresponds to the relation $\mathbf{A}_{,\mathbf{A}} = \mathbf{I} \boxtimes \mathbf{I}$ which, by virtue of definition (3.3), is equivalent to $\mathbf{C}^{IV}_{II}$. In the works of the second group the expression (5.17) corresponds to relation $\mathbf{A}^*_{,\mathbf{A}} = \mathbf{I} \otimes \mathbf{I}$ which, by virtue of definition (3.2), is equivalent to $\mathbf{C}^{IV}_{I}$ (note that in publication [13] the operation $\otimes$ is denoted as $\odot$). Taking into account that $\mathbf{C}^{IV}_{I} = (\mathbf{C}^{IV}_{II})^{ti}$ and $(\mathbf{C}^{IV}_{I})^{dr} = \mathbf{C}^{IV}_{I}$ as it follows from the schemes (1.1), expressions (3.2), (3.3) and schemes (5.2), we can conclude that the derivative of the tensor $\mathbf{A}$ with respect to the tensor $\mathbf{A}$ in the works of the second group coincides with an accuracy to the operations of transposition "$ti$" and "$dr$" with the similar derivative in the works of the first group. This suggests that the derivative $\mathbf{L}^{IV}_{II}$ is defined in terms of the derivative $\mathbf{L}^{IV}_{I}$ according to relations $(5.1)_1$, $(5.1)_2$ and $(5.1)_3$.

The derivative of the second rank tensor $\mathbf{A}^T$ with respect to the tensor $\mathbf{A}$ is defined by expression (5.18): $\mathbf{A}^T_{,\mathbf{A}} = \mathbf{C}^{IV}_{III}$. In the works of the first group, this expression corresponds to the relation $\mathbf{A}^T_{,\mathbf{A}} = \mathbb{T}$, where $\mathbb{T}$ is equivalent to $\mathbf{C}^{IV}_{III}$ in accordance with Table 1. In the works of the second group expression (5.18) corresponds to the relation $(\mathbf{A}^T)^*_{,\mathbf{A}} = \mathbb{T}$, but here according to Table 1 $\mathbb{T}$ is equivalent to $\mathbf{C}^{IV}_{II}$. Given that $[(\mathbf{C}^{IV}_{III})^{ti}]^{dr} = \mathbf{C}^{IV}_{II}$, we conclude that the derivative of the tensor $\mathbf{A}^T$ with respect to the tensor $\mathbf{A}$ in the works of the second group coincides with an accuracy to the operations of transposition "$ti$" and "$dr$" with the analogous derivative in the works of the first group, which suggests that the derivative $\mathbf{L}^{IV}_{II}$ should be defined in terms of derivative $\mathbf{L}^{IV}_{I}$ (see $(5.1)_1$, $(5.1)_2$ and $(5.1)_3$).

*6.5. Derivative* $\mathbf{A}^2_{,\mathbf{A}}$. The derivative $\mathbf{A}^2_{,\mathbf{A}} = \mathbf{C}^{IV}_{II} \overset{2}{*} \mathbf{A} + \mathbf{A} \cdot \mathbf{C}^{IV}_{II}$ is determined by relation (5.19). This derivative follows from (5.16) at $\mathbf{B} = \mathbf{S} \equiv \mathbf{A}$ taking into account (5.17). In the works of the first group it corresponds to the expression

$$\mathbf{A}^2_{,\mathbf{A}} = \mathbf{A} \boxtimes \mathbf{I} + \mathbf{I} \boxtimes \mathbf{A}^T, \qquad (6.1)_2$$

that follows from (6.1) (or $(6.1)_1$) at $\mathbf{B} = \mathbf{S} \equiv \mathbf{A}$. Since the equivalence of relations (5.16) and (6.1) (or $(6.1)_1$) has been demonstrated in subitem 6.3 it will suffice to show that $(6.1)_2$ indeed follows from $(6.1)_1$ at $\mathbf{B} = \mathbf{S} \equiv \mathbf{A}$, i.e. it is a special case of $(6.1)_1$. In view of the result obtained in subitem 6.4, $(6.1)_1$ can be written for the case under consideration in the following way:

$$\mathbf{A}^2_{,\mathbf{A}} = (\mathbf{A} \boxtimes \mathbf{I}) \cdot\cdot \mathbf{C}^{IV}_{II} \cdot\cdot \mathbf{C}^{IV}_{II} + (\mathbf{I} \boxtimes \mathbf{B}^T) \cdot\cdot \mathbf{C}^{IV}_{II} \cdot\cdot \mathbf{C}^{IV}_{II}. \qquad (6.1)_3$$

Taking into account that $\mathbf{C}^{IV}_{II} \cdot\cdot \mathbf{C}^{IV}_{II} = \mathbf{C}^{IV}_{III}$ (this can be readily shown, or none the worse, be found in [16], Annex I, §15), and also that $\mathbf{D}^{IV} \cdot\cdot \mathbf{C}^{IV}_{III} = \mathbf{D}^{IV}$ for any forth-rank tensor $\mathbf{D}^{IV}$ (the latter follows from the chain of equalities $\mathbf{D}^{IV} \cdot\cdot \mathbf{C}^{IV}_{III} = D^{ijkp}\mathbf{r}_i\mathbf{r}_j\mathbf{r}_k\mathbf{r}_p \cdot\cdot \mathbf{r}^m\mathbf{r}^n\mathbf{r}_n\mathbf{r}_m = D^{ijkp}\mathbf{r}_i\mathbf{r}_j\mathbf{r}_k\mathbf{r}_p = \mathbf{D}^{IV}$), we conclude that $(6.1)_3$, being a special case of $(6.1)_1$, reduces to $(6.1)_2$, which was what we set out to prove.

In the works of the second group relation (5.19) corresponds to the expression

$$(\mathbf{A}^2)^*_{,\mathbf{A}} = \mathbf{C}^{IV}_{I} \cdot \mathbf{A} + \mathbf{A} \cdot \mathbf{C}^{IV}_{I}, \qquad (6.2)_1$$

following from (6.2) at $\mathbf{B} = \mathbf{S} \equiv \mathbf{A}$ with account of the result of subpoint 6.4 ($\mathbf{A}^*_{,\mathbf{A}} = \mathbf{C}^{IV}_{I}$). The first term in $(6.2)_1$ is represented as $\mathbf{C}^{IV}_{I} \cdot \mathbf{A} = \mathbf{r}_m\mathbf{r}^m\mathbf{r}_n\mathbf{r}^n \cdot A^{ij}\mathbf{r}_i\mathbf{r}_j = \mathbf{r}_m\mathbf{r}^m A^{ij}\mathbf{r}_i\mathbf{r}_j = \mathbf{I} \otimes \mathbf{A}$. Applying successively to the next-to-last equality the transpose operations "$dr$" and "$ti$" (see schemes (5.2)), we obtain $[(A^{ij}\mathbf{r}_m\mathbf{r}^m\mathbf{r}_i\mathbf{r}_j)^{dr}]^{ti} = (A^{ij}\mathbf{r}_m\mathbf{r}^m\mathbf{r}_j\mathbf{r}_i)^{ti} = A^{ij}\mathbf{r}_m\mathbf{r}_j\mathbf{r}^m\mathbf{r}_i$, which, using the positional scalar multiplication operation introduced after relation (5.14), is transformed into $\mathbf{C}^{IV}_{II} \overset{2}{*} \mathbf{A}$. Indeed, $\mathbf{C}^{IV}_{II} \overset{2}{*} \mathbf{A} =$ $=\mathbf{r}_m\mathbf{r}^n\mathbf{r}^m\mathbf{r}_n \overset{2}{*} A^{ij}\mathbf{r}_i\mathbf{r}_j = \mathbf{r}_m(\mathbf{r}^n\cdot A^{ij}\mathbf{r}_i\mathbf{r}_j)\mathbf{r}^m\mathbf{r}_n = A^{ij}\mathbf{r}_m\mathbf{r}_j\mathbf{r}^m\mathbf{r}_i$. As a result we have $\mathbf{C}_I^{IV}\cdot\mathbf{A} = \mathbf{I}\otimes\mathbf{A} = [(\mathbf{C}_{II}^{IV} \overset{2}{*} \mathbf{A})^{ti}]^{dr}$. The second term in (6.2)$_1$ is represented as $\mathbf{A}\cdot\mathbf{C}_I^{IV} = A^{ij}\mathbf{r}_i\mathbf{r}_j\cdot\mathbf{r}^m\mathbf{r}_m\mathbf{r}^n\mathbf{r}_n = A^{ij}\mathbf{r}_i\mathbf{r}_j\mathbf{r}^n\mathbf{r}_n = \mathbf{A}\otimes\mathbf{I}$. The next-to-last equality can also be written as $[(\mathbf{A}\cdot\mathbf{C}_{II}^{IV})^{ti}]^{dr}$. Indeed, $\mathbf{A}\cdot\mathbf{C}_{II}^{IV} = A^{ij}\mathbf{r}_i\mathbf{r}_j\cdot\mathbf{r}^m\mathbf{r}_n\mathbf{r}_m\mathbf{r}^n = A^{ij}\mathbf{r}_i\mathbf{r}_n\mathbf{r}_j\mathbf{r}^n$. Taking into account the schemes of (5.2), $(\mathbf{A}\cdot\mathbf{C}_{II}^{IV})^{ti} = (A^{ij}\mathbf{r}_i\mathbf{r}_n\mathbf{r}_j\mathbf{r}^n)^{ti} = A^{ij}\mathbf{r}_i\mathbf{r}_j\mathbf{r}_n\mathbf{r}^n$ and $[(\mathbf{A}\cdot\mathbf{C}_{II}^{IV})^{ti}]^{dr} = (A^{ij}\mathbf{r}_i\mathbf{r}_j\mathbf{r}_n\mathbf{r}^n)^{dr} = A^{ij}\mathbf{r}_i\mathbf{r}_j\mathbf{r}^n\mathbf{r}_n$ and we obtain for (6.2)$_1$ the following equivalent representations: $(\mathbf{A}^2)^*_{,\mathbf{A}} = \mathbf{C}_I^{IV}\cdot\mathbf{A} + \mathbf{A}\cdot\mathbf{C}_I^{IV} = \mathbf{A}\otimes\mathbf{I} + \mathbf{I}\otimes\mathbf{A} = [(\mathbf{C}_{II}^{IV} \overset{2}{*} \mathbf{A} + \mathbf{A}\cdot\mathbf{C}_{II}^{IV})^{ti}]^{dr}$. Here, the last equality is written with consideration of the fact that the transposition of the sum is equal to the sum of the transposed quantities. This last equality satisfies the necessary requirement for the existence of relationship between the derivative $\mathbf{L}_{II}^{IV}$ and the derivative $\mathbf{L}_I^{IV}$ (see (5.1)$_1$, (5.1)$_2$ and (5.1)$_3$), by which the consistency of the derivatives $(\mathbf{A}^2)^*_{,\mathbf{A}}$ and $\mathbf{A}^2_{,\mathbf{A}}$ is satisfied.

*6.6. Derivative* $\mathbf{A}^{-1}_{,\mathbf{A}}$. The derivative $\mathbf{A}^{-1}_{,\mathbf{A}} = -\mathbf{A}^{-1}\cdot\mathbf{C}_{II}^{IV} \overset{2}{*} \mathbf{A}^{-1}$ is defined by expression (5.20). In the works of the first group this derivative corresponds to the expression $\mathbf{A}^{-1}_{,\mathbf{A}} = -\mathbf{A}^{-1}\boxtimes\mathbf{A}^{-T}$. Using for $\mathbf{A}^{-1}$ the coordinate representation $\mathbf{A}^{-1} = A^{ij}\mathbf{r}_i\mathbf{r}_j$, we obtain that

$$\mathbf{A}^{-1}\boxtimes\mathbf{A}^{-T} = A^{ij}A^{kp}\mathbf{r}_i\mathbf{r}_p\mathbf{r}_j\mathbf{r}_k, \tag{6.3}$$

which, as can be easily shown, coincides with the expression $\mathbf{A}^{-1}\cdot\mathbf{C}_{II}^{IV} \overset{2}{*} \mathbf{A}^{-1}$. In the works of the second group relation (5.20) corresponds to the expression $(\mathbf{A}^{-1})^*_{,\mathbf{A}} = -\mathbf{A}^{-1}\otimes\mathbf{A}^{-1}$, the right-hand part of which has (up to sign) coordinate representation $A^{ij}A^{kp}\mathbf{r}_i\mathbf{r}_j\mathbf{r}_k\mathbf{r}_p$. This tensor is obtained by a sequential transposition "$ti$" and "$dr$" of the right-hand part of equality (6.3), which is consistent with relation (5.20). As a result, $(\mathbf{A}^{-1})^*_{,\mathbf{A}} = -\mathbf{A}^{-1}\otimes\mathbf{A}^{-1} = -[(\mathbf{A}^{-1}\cdot\mathbf{C}_{II}^{IV} \overset{2}{*} \mathbf{A}^{-1})^{ti}]^{dr} = [(\mathbf{A}^{-1}_{,\mathbf{A}})^{ti}]^{dr}$ and the necessary requirement for the coincidence of derivative $\mathbf{L}_{II}^{IV}$ and $\mathbf{L}_I^{IV}$ is satisfied (see (5.1)$_1$, (5.1)$_2$ and (5.1)$_3$).

*6.7. Derivative* $(\boldsymbol{\Psi}(\mathbf{A})\,\boldsymbol{\Lambda}(\mathbf{A}))_{,\mathbf{A}}$. In our paper, this derivative is determined by relation (5.22): $(\boldsymbol{\Psi}(\mathbf{A})\,\boldsymbol{\Lambda}(\mathbf{A}))_{,\mathbf{A}} = \boldsymbol{\Lambda}(\mathbf{A})\,\boldsymbol{\Psi}(\mathbf{A})_{,\mathbf{A}} + \boldsymbol{\Psi}(\mathbf{A})\,\boldsymbol{\Lambda}(\mathbf{A})_{,\mathbf{A}}$. In the works of the first group this derivative corresponds to the expression $(\boldsymbol{\Psi}(\mathbf{A})\,\boldsymbol{\Lambda}(\mathbf{A}))_{,\mathbf{A}} = \boldsymbol{\Lambda}(\mathbf{A})\otimes\boldsymbol{\Psi}(\mathbf{A})_{,\mathbf{A}} + \boldsymbol{\Psi}(\mathbf{A})\,\boldsymbol{\Lambda}(\mathbf{A})_{,\mathbf{A}}$, which, when considering the remark that following equations (3.1), coincides with (5.22) completely. In the works of the second group $(\boldsymbol{\Psi}(\mathbf{A})\,\boldsymbol{\Lambda}(\mathbf{A}))^*_{,\mathbf{A}} = \boldsymbol{\Lambda}(\mathbf{A})\,\hat{\boxtimes}\,\boldsymbol{\Psi}(\mathbf{A})_{,\mathbf{A}} + \boldsymbol{\Psi}(\mathbf{A})\,\boldsymbol{\Lambda}(\mathbf{A})^*_{,\mathbf{A}}$. Taking into account the chain of equalities $[(\mathbf{A}\otimes\mathbf{B})^{ti}]^{dr} = (\mathbf{A}\boxtimes\mathbf{B})^{dr} = \mathbf{A}\,\hat{\boxtimes}\,\mathbf{B}$, where $\mathbf{A}$ and $\mathbf{B}$ are the second-rank tensors, we arrive at the equality $[((\boldsymbol{\Psi}(\mathbf{A})\,\boldsymbol{\Lambda}(\mathbf{A}))_{,\mathbf{A}})^{ti}]^{dr} = (\boldsymbol{\Psi}(\mathbf{A})\,\boldsymbol{\Lambda}(\mathbf{A}))^*_{,\mathbf{A}}$, which satisfies the requirement for the relationship between the derivative $\mathbf{L}_{II}^{IV}$ and the derivative $\mathbf{L}_I^{IV}$ ( see (5.1)$_1$, (5.1)$_2$ and (5.1)$_3$).

Table 3 summarizes the relations obtained in this article and in the publications by other authors, which have been discussed in the subitems of this Section.

**Table 3** *Correspondence between the rules for differentiation of the tensor argument functions obtained in different publications, and also between the resulting derivatives*

| Item | Proposed relations | Relations of the I group [8-10] | Relations of the II group [11-14] |
|---|---|---|---|
| 6.1. | $\Phi(\mathbf{A})_{,\mathbf{S}} = \Phi(\mathbf{A})_{,\mathbf{A}}\cdot\cdot\mathbf{C}_{II}^{IV}\cdot\cdot\mathbf{A}_{,\mathbf{S}}$ | $\Phi(\mathbf{A})_{,\mathbf{S}} = \Phi(\mathbf{A})_{,\mathbf{A}}(\cdot\cdot)\mathbf{A}_{,\mathbf{S}}$ | $\Phi(\mathbf{A})_{,\mathbf{S}} = \Phi(\mathbf{A})_{,\mathbf{A}}[\cdot\cdot]\mathbf{A}^*_{,\mathbf{S}}$ |
| 6.2. | $\boldsymbol{\Phi}(\mathbf{A})_{,\mathbf{S}} = \boldsymbol{\Phi}(\mathbf{A})_{,\mathbf{A}}\cdot\cdot\mathbf{C}_{II}^{IV}\cdot\cdot\mathbf{A}_{,\mathbf{S}}$ | $\boldsymbol{\Phi}(\mathbf{A})_{,\mathbf{S}} = \boldsymbol{\Phi}(\mathbf{A})_{,\mathbf{A}}(\cdot\cdot)\mathbf{A}_{,\mathbf{S}}$ | $\boldsymbol{\Phi}(\mathbf{A})^*_{,\mathbf{S}} = \boldsymbol{\Phi}(\mathbf{A})^*_{,\mathbf{A}}[\cdot\cdot]\mathbf{A}^*_{,\mathbf{S}}$ |
| 6.3. | $(\mathbf{A}(\mathbf{S}) \cdot \mathbf{B}(\mathbf{S}))_{,\mathbf{S}} =$ <br> $= \mathbf{A}_{,\mathbf{S}} \overset{2}{*} \mathbf{B} + \mathbf{A} \cdot \mathbf{B}_{,\mathbf{S}}$ | $(\mathbf{A}(\mathbf{S}) \cdot \mathbf{B}(\mathbf{S}))_{,\mathbf{S}} = (\mathbf{A} \boxtimes \mathbf{I})(\cdot\cdot)\mathbf{B}_{,\mathbf{S}} +$ <br> $+ (\mathbf{I} \boxtimes \mathbf{B}^T)(\cdot\cdot)\mathbf{A}_{,\mathbf{S}}$ | $(\mathbf{A}(\mathbf{S}) \cdot \mathbf{B}(\mathbf{S}))^*_{,\mathbf{S}} =$ <br> $= \mathbf{A}^*_{,\mathbf{S}} \cdot \mathbf{B} + \mathbf{A} \cdot \mathbf{B}^*_{,\mathbf{S}}$ |
| 6.4. | $\mathbf{A}_{,\mathbf{A}} = \mathbf{C}^{IV}_{II}$ <br> $\mathbf{A}^T_{,\mathbf{A}} = \mathbf{C}^{IV}_{III}$ | $\mathbf{A}_{,\mathbf{A}} = \mathbf{I} \boxtimes \mathbf{I}$ <br> $\mathbf{A}^T_{,\mathbf{A}} = \mathbf{I} \hat{\boxtimes} \mathbf{I}$ | $\mathbf{A}^*_{,\mathbf{A}} = \mathbf{I} \otimes \mathbf{I}$ <br> $\mathbf{A}^*_{,\mathbf{A}} = \mathbf{I} \boxtimes \mathbf{I}$ |
| 6.5. | $\mathbf{A}^2_{,\mathbf{A}} = \mathbf{C}^{IV}_{II} \overset{2}{*} \mathbf{A} + \mathbf{A} \cdot \mathbf{C}^{IV}_{II}$ | $\mathbf{A}^2_{,\mathbf{A}} = \mathbf{A} \boxtimes \mathbf{I} + \mathbf{I} \boxtimes \mathbf{A}^T$ | $(\mathbf{A}^2)^*_{,\mathbf{A}} = \mathbf{C}^{IV}_{I} \cdot \mathbf{A} + \mathbf{A} \cdot \mathbf{C}^{IV}_{I}$ |
| 6.6. | $\mathbf{A}^{-1}_{,\mathbf{A}} = -\mathbf{A}^{-1} \cdot \mathbf{C}^{IV}_{II} \overset{2}{*} \mathbf{A}^{-1}$ | $\mathbf{A}^{-1}_{,\mathbf{A}} = -\mathbf{A}^{-1} \boxtimes \mathbf{A}^{-T}$ | $(\mathbf{A}^{-1})^*_{,\mathbf{A}} = -\mathbf{A}^{-1} \otimes \mathbf{A}^{-1},$ |
| 6.7. | $(\boldsymbol{\Psi}(\mathbf{A})\boldsymbol{\Lambda}(\mathbf{A}))_{,\mathbf{A}} =$ <br> $= \boldsymbol{\Lambda}(\mathbf{A})\boldsymbol{\Psi}(\mathbf{A})_{,\mathbf{A}} +$ <br> $+ \boldsymbol{\Psi}(\mathbf{A})\boldsymbol{\Lambda}(\mathbf{A})_{,\mathbf{A}}$ | $(\boldsymbol{\Psi}(\mathbf{A})\boldsymbol{\Lambda}(\mathbf{A}))_{,\mathbf{A}} =$ <br> $= \boldsymbol{\Lambda}(\mathbf{A}) \otimes \boldsymbol{\Psi}(\mathbf{A})_{,\mathbf{A}} +$ <br> $+ \boldsymbol{\Psi}(\mathbf{A})\boldsymbol{\Lambda}(\mathbf{A})_{,\mathbf{A}}$ | $(\boldsymbol{\Psi}(\mathbf{A})\boldsymbol{\Lambda}(\mathbf{A}))^*_{,\mathbf{A}} =$ <br> $= \boldsymbol{\Lambda}(\mathbf{A}) \hat{\boxtimes} \boldsymbol{\Psi}(\mathbf{A})_{,\mathbf{A}} +$ <br> $+ \boldsymbol{\Psi}(\mathbf{A})\boldsymbol{\Lambda}(\mathbf{A})^*_{,\mathbf{A}}$ |

## VII. Conclusion

The rules for constructing derivatives of the scalar and tensor functions of the tensor argument with respect to a tensor and the form of representation of these derivatives depend on the accepted scheme of the interaction between the basis vectors during the operation of double scalar product of two tensors. Today, three groups of rules and representation forms have been developed and used. In the works of the first group, these rules are based on a widely used scheme of the basis vector interaction in a double scalar product of the second and fourth ranks tensors, which however, results sometimes in rather exotic forms of representation of derivatives, which essentially differ from representation of derivatives generally used in conventional mathematical analysis. In the works of the second group, the rules for constructing the derivatives of the scalar and tensor functions of tensor argument with respect to a tensor and the forms of their representation do not differ from the rules of constructing derivatives and the forms of their representation used in conventional mathematical analysis. However, in this case their construction is based on a special scheme of the interaction between the basis vectors in the double scalar product of the second- and fourth-rank tensor. The works of the third group also operate with a widespread scheme of the interaction between the basis vectors, which, however, is different from that used in the works of the first group, while the resulting rules for constructing derivatives of the scalar and tensor functions of the tensor argument with respect to a tensor and the form of derivative representation differ quite insignificantly from the rules for constructing derivatives and the forms of representation generally used in conventional mathematical analysis. For any of these groups, the derivatives of the scalar functions of the tensor argument with respect to a tensor are similar. The derivatives of the tensor functions of the tensor argument obtained in the works of the first and third groups represent the same tensor object which is the fourth-rank tensor $\mathbf{L}^{IV}_{I}$, determining the derivative $\boldsymbol{\Phi}(\mathbf{A})_{,\mathbf{A}}$, and differ only in the form of representation. A fourth-rank tensor $\mathbf{L}^{IV}_{II}$, which defines the derivative $\boldsymbol{\Phi}(\mathbf{A})^*_{,\mathbf{A}}$ of the tensor function of the tensor argument obtained in the works of the second group, differs from the tensor $\mathbf{L}^{IV}_{I}$ and is the result of successive transposition of "$ti$" and "$dr$" of the latter one. One should also keep in mind the peculiarities of the $[\cdot\cdot]$ operation, which is consistently used by the authors of the second group of publications. Taking all this into account, the ambiguity in the determination of derivatives of the tensor function of the tensor argument in the works of different groups, as noted earlier, will not be critical to analysis in solid mechanics, which are one of the main consumers of the results of this article.

## Acknowledgements

Author appreciate greatly assistance by L.V. Semukhina for preparation of English variant of this paper.